\documentclass{article}

\usepackage{shapepar,bm}
\usepackage[dvips,arrow,matrix,ps,color,line]{xy}
\usepackage{theorem,amsfonts,amssymb,amscd}
\usepackage{geometry}
\usepackage{makeidx}

\usepackage{graphicx}
\usepackage{graphics}

\usepackage{epstopdf}

\usepackage{amssymb, graphics,hyperref}

\hypersetup{
    colorlinks = true,
    linkcolor = blue,
    anchorcolor = red,
   citecolor = blue,
    filecolor = red,
    pagecolor = red,
    urlcolor = blue}

\newcommand{\mathsym}[1]{{}}

\usepackage[usenames]{color}
\definecolor{MyLightMagenta}{cmyk}{0.1,0.8,0,0.1}
\definecolor{MyDarkBlue}{rgb}{0.1,0,0.3}

\makeindex

\def\NN{\mathbb N}

\def\bfm{{\mathbf m}}

\def\bfb{{\mathbf b}}
\def\bfX{{\mathbf X}}

\def\Qcal{\mathcal Q}

\def\d{\partial}

\def\Gcal{\mathcal G}

\def\bfx{{\mathbf x}}

\def\ovDc{\overline{\Dcal}}

\def\ZZ{\mathbb Z}

\def\CC{\mathbb C}
\def\frakB{\mathfrak B}

\def\Ecal{{\mathcal E}}

\def\QQ{\mathbb Q}
\def\PP{\mathbb P}

\def\cocoa{{\hbox{\rm C\kern-.13em o\kern-.07em C\kern-.13em o\kern-.15em A}}}

\def\wb{{[\bfb]}}
\def\wbrl{b_{\lambda_r}\w b_{1+\lambda_{r-1}}\w\cdots\w b_{r-1+\lambda_1}}

\def\Dcal{\mathcal D}

\def\Res{{\rm Res}}

\def\Hom{{\rm Hom}}

\def\End{\mathrm{End}}

\def\sigm{{\sigma}}
\def\bsig{{\bm\sigma}}
\def\ovsig{{\overline{\sigma}}}
\def\blamb{{\bm \lambda}}
\def\bmu{{\bm\mu}}

\def\Pcal{{\mathcal P}}

\def\Vee{*}

\def\w2M{\bigwedge^2M}

\def\wM{\bigwedge M}

\def\w{\wedge }
\def\bw{\bigwedge }

\protect
\protect
\protect
\protect

\protect
\protect

\def\sra{\rightarrow}
\def\lra{\longrightarrow}

\def\proof{\noindent{\bf Proof.}\,\,}
\def\qedd{{\vrule height4pt width4pt depth0pt}\medskip}
\def\qed{{\hfill\vrule height4pt width4pt depth0pt}\medskip}

\def\be{\begin{equation}}
\def\ee{\end{equation}}
\def\bclm{\begin{claim}}
\def\eclm{\end{claim}}
\def\beqn{\begin{eqnarray}}
\def\eeqn{\end{eqnarray}}
\def\beqn*{\begin{eqnarray*}}
\def\eeqn*{\end{eqnarray*}}

\theoremstyle{change}
\theorembodyfont{\rmfamily}
\newtheorem{thm}{Theorem}[section]
\newtheorem{corol}{Corollary}[section]
\newtheorem{defin}{Definition}[section]

\newtheorem{prop}{Proposition}[section]
\newtheorem{lem}{Lemma}[section]
\newtheorem{rmk}{Remark}[section]
\newtheorem{claim}{}[section]

\global
\global

\def\no@breaks#1{{\def\\{ \ignorespaces}#1}}    

  \makeatletter
\def\cleardoublepage{\clearpage\if@twoside \ifodd\c@page\else
\hbox{} \thispagestyle{empty}
\newpage
\if@twocolumn\hbox{}\newpage\fi\fi\fi} \makeatother

\usepackage{eso-pic,graphicx}
\makeatletter
\newcommand\BackgroundPicture[2]{%
  \setlength{\unitlength}{1pt}%
  default \put(0,\strip@pt\paperheight){%
  \parbox[t][\paperheight]{\paperwidth}{%
    \vfill
     \centering \includegraphics[angle=#2, width=15cm, height=15cm,  bb=0 0 150 150]{#1}
    \vfill
}}} %
\makeatother
\date{}

\usepackage{eulervm}

\usepackage{amsmath}

\title{On  Pl\"ucker Equations Characterizing Grassmann Cones }

\author{\sc{Letterio Gatto, Parham Salehyan}}

\begin{document}

\maketitle
\date{}

\smallskip
\abstract {\noindent \,The KP hierarchy (after Kadomtsev and Petshiasvily) is a system of infinitely many PDEs in Lax form defining a universal family of iso-spectral deformation of an ordinary linear differential operator. It is a classical result due to Sato's japanese school that the rational solutions to the KP hierarchy are parametrized by a cone over an
 infinite-dimensional Grassmann variety. The present survey will revisit this fact from the point of view of {\em Schubert derivations} on a Grassmann algebra. These enable to encode the classical Pl\"ucker equations of  Grassmannians of $r$-dimensional subspaces in a formula whose limit
   for $r\sra \infty$ coincides with  the KP hierarchy, phrased in terms of  vertex operators,    showing in particular  how the latter is intimately related  to  Schubert calculus.}

\tableofcontents
\addcontentsline{toc}{section}{Introduction}

\section*{Introduction}

 \paragraph{} This survey article, which may also serve as  background
material while reading, for instance, \cite[Section 4]{Arb} and parts of~\cite{mulase,sato, segalw}, has the purpose to advertise the notion of {\em Schubert derivation} on an exterior algebra, introduced in~\cite{G1} (see also~\cite{GatSal2}), by showing  how it provides another approach to look at   the quadratic equations describing the Pl\"ucker embedding of  Grassmannians -- a very classical and widely studied subject. In particular,  it   allows i) to ``discover''  the vertex operators generating the fermionic vertex superalgebra (in the sense of~\cite[Section 5.3]{frbzvi}); ii) to compute their bosonic expressions as in~\cite[Lecture 5]{KR};  iii) to interpret them in terms of Schubert derivations and iv) to provide an almost effortless deduction of  the celebrated Hirota bilinear form of the KP hierarchy (after Kadomtsev and Petviashvilii) \cite{KR}.
Let $B:=\QQ[x_1,x_2,\ldots]$ and denote by $B_{(0)}$ its quotient field. A formal pseudo differential operator with coefficients in $B_{(0)}$ is a formal Laurent series $\sum_{i\leq n}a_{-i}(\bfx)\partial^{i}$, where $\d^{-1}$ denotes a formal inverse of the operator $\d:=\d/\d x_1$. The KP hierarchy concerns the evolution of  a first order normalized pseudo differential operator, $L:=\d+u_1(\bfx)\d^{-1}+\cdots\in B_{(0)}((\d^{-1}))$ obeying  the Lax equations 
\be
{d L\over d x_n}=[(L^n)_+,L],\label{eq:lax}
\ee
where $(L^n)_+$ denotes the differential part of the $n$-th power of $L$ -- a linear ordinary differential operator of order $n$. 
Equations~(\ref{eq:lax}) arise as compatibility conditions for an isospectral deformation of an ordinary differential operator \cite[Section 4]{mulase}. It is easily seen that the first two non trivial equations of the hierarchy ($n=2,3$) yield the celebrated KP equation
\[
{3\over 4}f_{yy}-\left(f_t-{1\over 4}f_{xxx}-3ff_{x}\right)_x=0,
\]
upon identifying $x_1=x$, $x_2=y$, $x_3=t$ and $f:=u_1$ \cite[Section 4]{mulase}.
It is  a fundamental observation due to Sato~\cite{sato,sato1}, and widely developed by his Kyoto school \cite{DJKM1,DJKM2,JM},  that solutions $L$ to the KP hierarchy~(\ref{eq:lax}) are parameterized by the points of a Grassmannian parametrizing infinite dimensional subspaces of an infinite dimensional vector space. Indeed, 
explicit solutions to the KP equation in Lax form can be constructed by {\em tau} functions  which, for the limited purposes of this paper, we define just  as those  polynomials $\tau\in B$ satisfying
\be
\Res_z\big(\Gamma^\Vee(z)\tau\otimes\Gamma(z)\tau\big)=0,\label{eq:KPHwVOx}
\ee
where $z$ is a formal variable, $\Res_z$ denotes the coefficient of $z^{-1}$ of a formal Laurent series and  $\Gamma(z),\Gamma^\Vee(z):B\sra B((z))$ are the ``bosonic vertex operators'':
\be
\Gamma(z):=\exp(\sum_{i\geq 1}x_iz^i)\exp\left(\hskip-2pt -\sum_{i\geq 1}{1\over iz^i}{\d\over \d x_i}\right)\,\,\, \mathrm{and}\,\,\, \Gamma^\Vee(z):=\exp(-\sum_{i\geq 1}x_iz^i)\exp\left(\hskip-2pt\sum_{i\geq 1}{1\over iz^i}{\d\over \d x_i}\right)\hskip-1pt.\label{eq:inteqvo}
\ee
For a {\em tau} function $\tau\in B$, 
let
\[
P_\tau(z)=\tau(\bfx)^{-1}\exp\left(-\sum_{i\geq 1}{1\over iz^i}{\d\over \d x_i}\right)\tau(\bfx)\in B_{(0)}[z^{-1}].
\]
Define  $P_\tau(\d)$ to be the evaluation of $P_\tau(z)$ at $z^{-1}=\d^{-1}$. Then $L_\tau=P_\tau(\d)^{-1}\cdot \d\cdot P_\tau(\d)$ is a normalized formal pseudo--differential operator that satisfies the Lax equations~(\ref{eq:lax}) \cite[Section 4]{Arb}.

The point we shall focus on is that~(\ref{eq:KPHwVOx}), also known as {\em Hirota bilinear form of the KP hierarchy}, and which we refer to as the KP hierarchy {\em tout court},  encodes the Pl\"ucker equations of the cone  of decomposable tensors of  a semi-infinite exterior power of an infinite--dimensional vector space. 
This fact is mentioned and/or explained in a number of different ways, e.g.  in~\cite{alekalatsu, Arb, enolHarn,glimpses2,KL,mulase,mulase1,loopgroups,segalw} and surely in many more references. From our part we shall recover   expression~(\ref{eq:KPHwVOx})  by considering the  limit for $r\sra \infty$ of the main formula we prove in the present article, that  characterizes  Grassmann cones  of decomposable tensors in an $r$-th exterior power ($r\geq 1$).

 \paragraph{} For this purpose,  let $M_0$ be a free abelian group of infinite countable rank with a basis $\frakB_0:=(b_0,b_1,\ldots)$ and let $(\beta_j)_{j\geq 0}$ be the basis of the restricted dual $M_0^\Vee$ such that $\beta_j(b_i)=\delta_{ij}$. Let $\bw^rM_0$ be the $r$-th exterior power of $M_0$.
{\em The Grassmann cone}   $\Gcal_r$ of $\bw^rM_0$  is  the image of the (non-surjective) multilinear alternating map $M_0^r\sra \bw^rM_0$, given by    $(m_1,\ldots, m_r)\mapsto m_1\w\cdots\w m_r$.

\paragraph{}
It is well-known, see e.g.~\cite[Section 4]{Arb},   that $\bfm\in \bw^rM_0$ belongs to $\Gcal_r$ if and only if
\be
\sum_{i\geq 0}(\beta_i\lrcorner \bfm)\otimes (b_i\w \bfm)=0, \label{eq:clasdec}
\ee
where $\beta_i\lrcorner\bfm\in \bw^{r-1}M_0$ denotes the contraction of $\bfm$ against $\beta_i$ (Section~\ref{sccontra}). Equation~(\ref{eq:clasdec}) is equivalent to
\be
\Res_{z}\Big(\sum_{i\geq 0}(\beta_iz^{-i-1}\lrcorner\bfm)\otimes \sum_{j\geq 0}(b_jz^j\w\bfm)\Big)=0,\label{eq:clasdecres}
\ee
a trick we learned in ~\cite[Section 7.3]{KR}. We combine~(\ref{eq:clasdecres}) with the following observation: there are unique formal power series
 $\sigm_{\pm}(z):=\sum_{i\geq 0}\sigm_{\pm i}z^{\pm i}\in \End_\ZZ(\wM_0)[[z^{\pm 1}]]$ such that 
 \[
\left\{\begin{matrix}\sigm_{\pm}(z)(\bfm_1\w\bfm_2)&=&\sigm_{\pm}(z)\bfm_1\w\sigm_{\pm}(z)\bfm_2,& (\forall \bfm_1,\bfm_2\in \wM_0)&\cr\cr
\sigm_{\pm i}b_j&=& b_{j\pm i},\,\,\,\,\,\,\,\,\,\,\,\,\,\,\,\,\,\,\,\,\,\,\,\,\,\,\,\,\,\,\,\,\,\,\,\,\,\,\,\,\,\,\,\,\,\end{matrix}\right.
\]
by agreeing that $b_k=0$ if $k<0$.  Let $\ovsig_\pm(z):=\sum_{i\geq 0}\ovsig_{\pm i}z^{\pm i}$
 be the inverse in $\End_\ZZ(\wM_0)[[z^{\pm 1}]]$  of $\sigm_\pm(z)$.

\smallskip
The  main result of this paper is:
 \begin{thm}\label{thm01} {\em An element  $\bfm\in\bw^rM_0$ belongs to the Grassmann cone $\Gcal_r$ if and only if the  equality
\begin{center}
\begin{tabular}{|c|}\hline\\
\,\,\,\,\,\,\,\,\,\,\,\,$\Res_{z}\, \big[\ovsig_+(z)\ovsig_{-r+1}(\beta_0\lrcorner \sigm_-(z)\bfm)\otimes \sigm_+(z)\ovsig_{-}(z)(b_0\w \ovsig_{r}\bfm)\big]=0$\,\,\,\,\,\,\,\,\,\,\,\,\\ \\
\hline
\end{tabular}
\end{center}
 holds in $\bw^{r-1}M_0\otimes_\ZZ\bw^{r+1}M_0$.
}
\end{thm}

\smallskip
{Following~\cite{G1}, we have called {\em Schubert derivations} the maps $\sigm_{\pm }(z)$ as well as their formal multiplicative inverses in $\End_\ZZ(\bw M_0)$. The reasons for such terminology are that i) they are derivations of $\wM_0$ in the sense of Hasse and Schmidt (see~\cite{G1,GatSal2}) and ii) they  satisfy suitable Pieri and Giambelli formulas (Section~\ref{secshur1}). It follows, in particular, that the exterior power $\bw^r\CC^n$ can be regarded as  an irreducible module over the cohomology ring $H^*(G(r,n),\ZZ)$ of the Grassmannian  of $r$-planes in
    $\CC^n$ (\cite{G1,G2}), whence  the identification $H^*(G(r,n),\ZZ)\cong \bw^rH^*(\PP^{n-1})$,  also known in the literature as  {\em Satake isomorphism} \cite{GALGOLHI,GOLMA}.
    The reason why  Theorem~\ref{thm01} is interesting is that its proof is cheap and easily returns the equation of the KP hierarchy by letting  $r$ tending to $\infty$. }

\paragraph{} To state a relevant consequence of Theorem~\ref{thm01} above, we need to introduce a few more  pieces of notation.
  Let  $\Pcal_r:=\{\blamb:=(\lambda_1,\ldots,\lambda_r)\in\NN^r\,|\, \lambda_1\geq\cdots\geq \lambda_r\}
$ be the set of all {\em partitions of length at most $r$} and  $\Pcal:=\cup_{r\geq 0}\Pcal_r$  the set of all partitions.  For $\blamb\in\Pcal_r$,  let $\wb^r_\blamb:=b_{\lambda_r}\w b_{1+\lambda_{r-1}}\w\cdots\w b_{r-1+\lambda_1}$ so that $\bw^r\frakB_0:=(\wb^r_\blamb\,\big|\, \blamb\in\Pcal_r)$ is a $\ZZ$-basis of  $\bw^rM_0$. Then each $\bfm\in\bw^rM_0$
   can be uniquely written as a finite linear combination of the form $\sum_{\blamb\in\Pcal_r}a_\blamb\wb^r_\blamb$.
If $B_r:=\ZZ[e_1,\ldots, e_r]$ is the polynomial ring in the $r$ indeterminates $(e_1,\ldots, e_r)$, let
  $E_r(z):=1-e_1z+\cdots+(-1)^re_rz^r\in B_r[z]$. Construct a sequence  $H_r:=(h_j)_{j\in\ZZ}$ of elements of $B_r$ via the equality
\[
\sum_{i\in\ZZ}h_iz^i:=\sum_{n\geq 0}(1-E_r(z))^n,
\]
understood in  the ring of formal Laurent series $B_r[[z^{-1},z]]$. By construction,  $h_j=0$ if $j<0$, $h_0=1$ and, for $i>0$, $h_i$
 is a homogeneous polynomial of degree $i$ in $(e_1,\ldots, e_r)$ provided that, for all  $1\leq j\leq r$, each $e_j$ is given weight $j$.
   The Schur determinant associated to the sequence $H_r$ and to $\blamb\in\Pcal_r$ is by definition
\[
\Delta_\blamb(H_r):=\det(h_{\lambda_j-j+i})_{1\leq i,j\leq r}\in B_r.
\]
Using the  well--known fact  that  $B_r=\bigoplus_{\blamb\in\Pcal_r}\ZZ\cdot \Delta_\blamb(H_r)$, the map  
\be
\phi_r:B_r\sra \bw^rM_0,\label{eq0:mapfr}
\ee
  given by $\Delta_\blamb(H_r)\mapsto \wb^r_\blamb$, defines an  isomorphism of abelian groups (the ``boson-fermion correspondence of order  $r$'').  It enables to equip $\bw^rM_0$ with a structure of free $B_r$-module
of rank $1$, generated by $\wb^r_0$, that we shall denote by $\bw^rM_r$. We  regard  $\sigm_-(z),\ovsig_-(z)$
as  maps $B_r\sra B_r[z^{-1}]$ as well,   by defining   $\sigm_-(z)\Delta_\blamb(H_r)=\phi_r^{-1}(\sigm_-(z)\wb^r_\blamb)$
and $\ovsig_-(z)\Delta_\blamb(H_r)=\phi_r^{-1}(\ovsig_-(z)\wb^r_\blamb)$.  Denote by  $\sigm_-(z)H_r$ and $\ovsig_-(z)H_r$,  respectively,
 the sequences $(\sigm_-(z)h_j)_{j\in\ZZ}$ and $(\ovsig_-(z)h_j)_{j\in\ZZ}$ in $B_r[z^{-1}]$. 
The following equalities (Proposition~\ref{prophnz})
   \[
\sigm_-(z)h_n=\sum_{j=0}^nh_{n-j}z^{-j}\qquad  \mathrm{and} \qquad \ovsig_-(z)h_n=h_n-h_{n-1}z^{-1}
\]
hold in $B_r[z^{-1}]$,  for all $r\geq 1$, and Theorem~\ref{thm01} admits the following rephrasing:

\begin{thm}\label{thm02} {\em The element $\bfm:=\sum_{\blamb\in\Pcal_r}a_\blamb\wb^r_\blamb\in\bw^rM_0$
  belongs to $\Gcal_r$ if and only if the  equality below holds in $B_{r-1}\otimes_\ZZ B_{r+1}$
 \be
 \begin{tabular}{|c|}\hline\\
$\Res_z \sum_{\blamb,\bmu\in\Pcal_r}a_\blamb
a_\bmu E_{r-1}(z)\Delta_\blamb(\sigm_-(z)H_{r-1})\otimes
\displaystyle{1\over E_{r+1}(z)}\ovsig_-(z)\Delta_\bmu(H_{r+1})=0$.\\ \\ \hline
 \end{tabular}
  \label{eq1:mnthm}
  \ee
}

\end{thm}

\smallskip
\noindent
For example, for $r=2$ one may easily recover the equation of the Klein's quadric cutting out the Grassmannian $G(2, 4)$ in its Pl\"ucker embedding (Section~\ref{SecEx}). Other examples are discussed  in  the book~\cite{GatSal2}. They all  indicate that even for detecting the  Grassmann cone $\Gcal_2$, computations are quite painful, surely not as easy as  checking the simpler condition $\bfm\w \bfm=0$.
What makes Theorem~\ref{thm02} interesting, however, is that, on one hand,  the maps $\Gamma_r(z):B_r\sra B_{r+1}((z))$ and $\Gamma^\Vee_r(z):B_r\sra B_{r-1}((z))$, defined by
\[
\Gamma_r(z)\Delta_\blamb(H_r):={1\over E_{r+1}(z)}\ovsig_-(z)\Delta_\blamb(H_{r+1})\quad \mathrm{and}\quad \Gamma^\Vee_r(z)\Delta_\blamb(H_r):=E_{r-1}(z)\Delta_\blamb(\sigm_-(z)H_{r-1})
\]
and occurring in formula~(\ref{eq1:mnthm}),
 are precisely truncated versions of the vertex operators displayed in~(\ref{eq:inteqvo}) and, on the other hand, the Schubert derivations $\sigm_-(z)$ and $\ovsig_-(z)$ are well defined also for $r=\infty$. 
 
More precisely, let $(e_1,e_2,\ldots)$ be a sequence of infinitely many indeterminates and
 $E_\infty(z):=\sum_{i\geq 0}(-1)^ie_iz^i$. For all  $p\in B_\infty$ there is $r\geq 0$ such that $p\in B_s$ for all $s\geq r$. 
  We have:

\begin{thm}\label{mncor01} {\em An element  $\tau\in B_\infty:=\ZZ[e_1,e_2,\ldots]$  is solution to the equation
 \be
\Res_{z}\, \left( E_{\infty}(z)\sigm_-(z)p\otimes_\ZZ{1\over E_{\infty}(z)}\ovsig_-(z)p\right)=0 \label{eq1:mnthmm}
 \ee
 if and only if there exists $s\geq 0$ such that $\phi_r(\tau)\in\Gcal_r$  for all $r\geq s$. 
}
\end{thm}

\medskip
\noindent
It turns out that   equation~(\ref{eq1:mnthmm}) expresses the KP hierarchy~(\ref{eq:KPHwVOx}) over the integers. It has been  obtained by  using the  indeterminates $e_i$ and $h_j$ (which may be interpreted as elementary and complete symmetric polynomials), often more convenient than the variables $(x_i)$ used in~(\ref{eq:KPHwVOx}) and~(\ref{eq:inteqvo}). This is also sanctioned in the couple of important and  relatively recent articles~\cite{frpenkser,JingRozh}.




\paragraph{} By abuse of notation, let us write $\Gcal_r\otimes\QQ$ for the Grassmann cone of decomposable tensors of $\bw^r(M_0\otimes_\ZZ\QQ)$.
Let  $(x_1,x_2,\ldots)$ be the sequence of indeterminates over $\QQ$, implicitly defined by the equality:
 \[
 \exp(-\sum_{i\geq 1}x_iz^i):=E_\infty(z).
 \]
An immediate check shows that 
$B:=B_\infty\otimes \QQ=\QQ[x_1,x_2,\ldots]$. 


\begin{corol}\label{mncor02} {\em An element $\tau\in B$ is a solution to the KP hierarchy~{\em (\ref{eq:KPHwVOx})} if and only if there exists  $s\geq 0$ such that  $\phi_r(\tau)\in\Gcal_r\otimes\QQ$ for all $r\geq s$.
}

\end{corol}
\paragraph{} The paper is organized as follows. Section~\ref{preli} sets the preliminaries and notation used throughout the paper.
 Section~\ref{sechsder} recalls a few facts concerning Hasse-Schmidt ($HS$) derivations on exterior algebras as introduced in~\cite{G1} and treated in more details
  in~\cite{GatSal2}. The section proclaims the most powerful tool of the theory which we call (as in~\cite{G1}) {\em integration by parts formula}. That  the transpose
     of an $HS$-derivation is an  $HS$-derivation as well is also shown in Section~\ref{sechsder},  a fact heavily used   to prove Theorem~\ref{thm01}.
 Schubert derivations are studied in Section~\ref{secSchuder}, where  a few technical lemmas leading to
    the approximation  $B_r\sra B_r((z))$ of the vertex operators  are discussed.      
       A pivotal aspect of Section~\ref{GCPR} is that the Schubert derivations $\sigm_-(z)$ and its inverse $\ovsig_-(z)$ enjoy a stability property enabling
        to  define them as maps $B_r\sra B_r[z^{-1}]$. Their limit for $r\sra \infty$ give rise to the  ring homomorphisms $B\sra B[z^{-1}]$, which enter in the expression of the vertex operators.
        The  crucial property that  $\sigm_{- }(z),\ovsig_-(z)$ commute with taking Schur determinants, proven in Section~\ref{GCPR},  is obtained by exploiting  a powerful
         determinantal formula  due to Laksov and Thorup~\cite{LakTh1}. Section~\ref{SecEx} is entirely devoted to the standard example of decomposable tensors in a second wedge power,  faced  via  Theorem~\ref{thm02}. Eventually, Section~\ref{SecKP} is concerned with the limit of formula~(\ref{eq1:mnthm}) for $r\sra \infty$, whereby  $\bw^\infty M_0$ shall be understood as the projective limit $\underset{\leftarrow}{\lim}\bw^rM_0$ in the category of graded modules.
There,  instead of (re)showing that {\em tau}-functions   correspond to decomposable tensors in some infinite exterior power,  we  use  them  to define the Grassman cone $\Gcal_\infty$ as the locus of $\bfm\in \bw^\infty M_0$ such that $\phi_\infty^{-1}(\bfm)$ solves the $KP$ hierarchy, where $\phi_\infty$ is the analogous of~(\ref{eq0:mapfr}) for $r=\infty$.  \vspace{-1pt}
\bigskip

We are grateful to S.~G.~Chiossi, A.~Kasman, A.~Ricolfi and I.~Scherbak for useful comments and corrections and to C.~Araujo,
 M.~Jardim, P.~Mulassano and P.~Piccione for unreserved support.  Thanks are also due to the anonymous  referee for his/her patient careful reading and remarks,   that helped us to improve the shape of the article. 
 
 \smallskip
 {\em This paper is dedicated to Professor Piotr Pragacz, a living example of mathematical experience, on the occasion of his sixtieth birthday.}

 \section{Preliminaries and Notation}\label{preli}
This section is to fix notation and to list the pre-requisites we shall need in the sequel.
 \claim{}  We denote by $\Pcal$ the set of all  monotonic non increasing sequence $\blamb:=(\lambda_1\geq \lambda_2\geq\cdots)$ of non negative integers
  all zero but finitely many, said {\em partitions}. The {\em length} $\ell(\blamb)\in \Pcal$   is  $\sharp\{i\,|\, \lambda_i\neq 0\}$, the number
   of its non-zero {\em parts};  its {\em weight}  $|\blamb|:=\sum_{i\geq 1}\lambda_i$. 
   We shall denote by $\Pcal_r$ the set of all
    partitions of length at most $r$ and  by 
    \[
    \Pcal_{r,n}:=\{\blamb:=(\lambda_1,\ldots, \lambda_r)\in\Pcal_r\,|\, \lambda_1\leq n-r\}
    \]
     the elements of $\Pcal_r$ bounded by $n-r$. An expression $(1^{i_1},2^{i_2},\ldots, )$ such that all $i_j$ are zero, but finitely many, denotes the partition having $i_j$ parts equal to $j$.  The {\em Young
diagramme}\index{Young diagramme} of
$\blamb:=(\lambda_1,\ldots,\lambda_r)\in \Pcal$ is an array $Y(\blamb)$ of left-justified rows of boxes, such that the {\em j}-th row has $\lambda_j$ boxes, for $1\leq j\leq r$. 
The set $\Pcal_r$ is a monoid with respect to the sum
      $\blamb+\bmu=(\lambda_i+\mu_i)_{i\geq 1}$ whose neutral element is the  null partition $(0)$ with all the parts equal
       to zero. If $\blamb,\bmu\in\Pcal_r$,  we shall write $\bmu\subseteq \blamb$ if the Young diagram of $\bmu$ is contained
        in the Young diagram of $\blamb$.

 \claim{}
We shall denote by $M_0$ the free abelian group $\ZZ[X]$ and by $\frakB_0:=(b_0,b_1,\ldots)$ its standard basis
 $(1,X,X^2,\ldots)$. Let $\wM_0:=\bigoplus_{r\geq 0}\bw^rM_0$ be the exterior algebra of $M_0$. Then $\bw^0M_0:=\ZZ$ and,
  for all $r\geq 1$,  the $r$-th exterior power of $M_0$ is the $\ZZ$-linear span of $\bw^r\frakB_0:=(\wb^r_\blamb\,|\, \blamb\in\Pcal_r)$,
   where  $\wb^r_\blamb:=\wbrl$. In particular $\wb^r_0:=\wb^r_{(0)}=b_0\w b_1\w\cdots\w b_{r-1}$.

\claim{}\label{sccontra} Let $\bfm\in\bw^rM_0$, $r\geq 1$. Its contraction, $\beta\lrcorner\bfm$, against $\beta\in M_0^\Vee$ is the
 unique element of $\bw^{r-1}M_0$ such that
\be
\gamma(\beta\lrcorner\bfm)=(\beta\w \gamma)(\bfm),\qquad \forall \gamma\in\bw^{r-1}M_0^\Vee. \label{eq4:contra}
\ee
It turns out that $\beta\lrcorner: (\wM_0,\wedge)\sra (\wM_0,\wedge)$ is the unique derivation of degree $-1$ such that $\beta\lrcorner m=\beta(m)$,
  for all $m\in M_0$.

\claim{}\em Let $B_0=\ZZ$ and for $r\geq 1$ denote by $B_r$ the polynomial
ring $\ZZ[e_1,\ldots,e_r]$. Accordingly, we let $E_0(z)=1$ and
$E_r(z)=1-e_1z+\cdots+(-1)^re_rz^r\in B_r[z]$  for $r\geq 1$.  The equality
\[
\sum_{n\in\ZZ}h_nz^n={1\over E_r(z)}=\sum_{i\geq 0}(1-E_r(z))^i,
\]
read in the abelian group of formal Laurent series $B_r[[z^{-1},z]]$, defines the bilateral sequence $H_r:=(h_j)_{j\in\ZZ}$ of elements of $B_r$.

Then,  by construction,  $h_j=0$ if $j<0$, $h_0=1$ and for all $j>0$, $h_j$ is a polynomial in $e_1,\ldots, e_r$ of weighted degree $j$,
  after declaring that $e_i$ is given degree $i$.
\begin{rmk}\label{2notrmkhrn}  For
the terms of the sequence $H_r$,
the more careful notation $h_{r,n}$ should be preferred, in place of just $h_n$, to keep track of their dependence on $r$. To make the notation less heavy we decided however to drop the subscript {r}, hoping for  the context
being sufficient to avoid confusions. 
\end{rmk}

\claim{}\label{prelwtgrad} It is well known  that (e.g.~\cite[p.~41]{MacDonald})
$
B_r=\bigoplus_{\blamb\in\Pcal_r}\ZZ\cdot \Delta_\blamb(H_r),
$
where
\[
\Delta_\blamb(H_r):=\det(h_{\lambda_j-j+i})_{1\leq i,j\leq r}.
\]
  Each partition
   $\blamb':=(1^{i_1}\ldots r^{i_r})$  is {\em conjugated} to  a partition $\blamb\in\Pcal_r$, of weight $i_1+2i_2+\cdots+ri_r$,
    whose Young diagramme is the transpose of $Y(\blamb')$. Define $e^\blamb:=e_1^{i_1}\cdots e_r^{i_r}$ and let $(B_r)_w:=\bigoplus_{|\blamb|=w}\ZZ\cdot e^\blamb$.
The direct sum decomposition $B_r:=\bigoplus_{w\geq 0}(B_r)_w$,  shall be referred to as the {\em weight graduation} of $B_r$. Then $(B_r)_w$ is the submodule of $B_r$ of all the polynomials of weight $w$ and it coincides with $\bigoplus_{|\blamb|=w}\ZZ\cdot \Delta_\blamb(H_r)$.

\vspace{-1pt}
 \section{Hasse-Schmidt Derivations on an Exterior Algebra} \label{sechsder}

\claim{} Given any module $M$ over a commutative ring $A$ with unit, there is  an obvious $A$-module isomorphism
\be
\End_A(M)[[z]]\sra{\mathrm{Hom}}_A(M, M[[z]]),\label{eq3:iso}
\ee
 where $\End_A(M)[[z]]$  denotes the formal power series with $\End_A(M)$-coefficients  in an indeterminate $z$.
  If $\Dcal(z)\in\End_A(M)[[z]]$,  we denote in the same way its image through the map~(\ref{eq3:iso}).
\begin{defin}{\em A Hasse--Schmidt {\em ($HS$)} derivation of $\wM$ is an algebra homomorphism
 \linebreak  $\Dcal(z):\wM\sra\wM[[z]]$, i.e:
\be
\Dcal(z)(\bfm_1\w\bfm_2)=\Dcal(z)\bfm_1\w\Dcal(z) \bfm_2, \quad \forall \bfm_1,\bfm_2\in\wM.\label{eq2:HSDA}
\ee
}
\end{defin}
Define a  sequence $\Dcal:=(D_0,D_1,\ldots)$ of endomorhisms of $\wM$ through the equality:
\[
\sum_{j\geq 0}D_j\bfm\cdot z^j:=\Dcal(z)\bfm.
\]
Then equation~(\ref{eq2:HSDA}) holds if and only if  the sequence $\Dcal$ obeys the {\em higher order Leibniz rules}:
\be
D_i(\bfm_1\w\bfm_2)=\sum_{j=0}^iD_i\bfm_1\w D_{i-j}\bfm_2,\qquad i\geq 0.
\ee
In particular $D_0$ is an algebra homomorphism, and   $D_1$ is a (usual) derivation of the $\w$-algebra $\wM_0$, provided that $D_0={\mathrm id}_{\wM}$. Moreover, if $D_0$ is an
 automorphism of $\wM$, it turns out that  the formal power series $\Dcal(z)$ is invertible in $\End_A(\wM)[[z]]$. Denote by $\ovDc(z)$ its inverse.

\begin{prop}[see \cite{G1}]\label{propinv} {\em The set $HS(\wM)$ of all $HS$-derivations on $\wM$  is a subalgebra of $\End_\ZZ(\wM)[[z]]$, with respect to the product
\[
\Dcal(z)\Ecal(z)=\sum_{j\geq 0}\sum_{i=0}^j (D_i\circ E_{j-i})z^{j},
\]
where $\sum_{i\geq 0}D_iz^i:=\Dcal(z)$ and $\sum_{i\geq 0}E_iz^i:=\Ecal(z)$.
In particular
   if $D_0$ is an automorphism of $\wM$, then the inverse formal power series $\ovDc(z)$  is an $HS$-derivation if and only if $\Dcal(z)$ is. }
\end{prop}

\begin{prop} \label{unicitypr} {\em Let $f(z)\in \End_A(M)[[z]]$. There exists a unique $HS$-derivation $\Dcal^f(z):\wM\sra \wM[[z]]$ such that
\[
\Dcal^f(z)m=f(z)(m)
\]
for all $m\in M$.
}
\end{prop}

\proof Let us first prove the existence. For all $r\geq 1$, consider the unique $A$-linear extension of the  map
$\Dcal^f(z):M^{\otimes r}\sra \bw^rM[[z]]$  defined by
$
\widehat{f}(z)(m_1\otimes\cdots\otimes m_r)=f(z)m_1\w\cdots\w f(z)m_r.
$
This map is clearly alternating and thus it factorizes through a
unique homomorphism $\Dcal^f_r(z):\bw^rM\sra\bw^rM[[z]]$ such that
$\Dcal^f_r(z)(m_1\otimes\cdots\otimes m_r)=f(z)m_1\w\cdots\w
f(z)m_r$. Each $\bfm\in\wM$ is a finite sum
$\bfm_{1}+\cdots+\bfm_{s}$ of homogeneous elements, i.e.
$\bfm_i\in\bw^iM$ for some $i\geq 0$ (notice that $\Dcal^f(z)a=a$ for all $a\in A$). Define $\Dcal^f(z)\bfm$ as
$\sum_{i=1}^s\Dcal^f_i(z)\bfm_i$. We want to show that for
$\bfm_1,\bfm_2\in \wM$
\[
\Dcal^f(z)(\bfm_1\w\bfm_2)=\Dcal^f(z)\bfm_1\w \Dcal^f(z)\bfm_2.
\]
Without loss of generality, we may assume they are homogeneous with respect to the graduation of $\wM$, i.e.
 $\bfm_1=m_{11}\w\cdots\w m_{1r}$ and $\bfm_2:=m_{21}\w\cdots\w m_{2s}$. Thus
\begin{eqnarray*}
\Dcal^f(z)(\bfm_1\w\bfm_2)&=&\Dcal^f(z)(m_{11}\w\cdots\w m_{1r}\w m_{21}\w\cdots\w m_{2s})\\
&=&f(z)m_{11}\w\cdots\w f(z)m_{1r}\w f(z)m_{21}\w\cdots\w f(z)m_{2s}=\\
&=&(f(z)m_{11}\w\cdots\w f(z)m_{1r})\w (f(z)m_{21}\w\cdots\w f(z)m_{2s})\\
&=&\Dcal^f(z)\bfm_1\w\Dcal^f(z)\bfm_2.
\end{eqnarray*}
To prove unicity, let $\widehat{D}(z)$ be any $HS$-derivation on $\wM$ such that $\widehat{D}(z)m=f(z)m$ for all $m\in\ M$.
 Then for all homogeneous element $\bfm:=m_1\w\cdots \w m_r\in\bw^rM$:
\begin{eqnarray*}
\hskip46pt \widehat{D}(z)\bfm&=&\widehat{D}(z)(m_1\w\cdots \w m_r)=\widehat{D}(z)m_1\w\cdots\w\widehat{D}(z)m_r\\
&=&f(z)m_1\w\cdots\w f(z)m_r=\Dcal^f(z)(m_1\w\cdots\w m_r)=D^f(z)\bfm.\hskip46pt \qed
\end{eqnarray*}

The main tool of the paper is the following observation for which we omit the totally obvious proof. It  is responsible, in our context, of the emergence of the vertex operators.
\begin{prop}[{\bf Integration by Parts}] {\em Assume that  $\Dcal(z)\in HS(\wM)$ is invertible in the sense of Proposition~{\em \ref{propinv}}. Then  the {\em integration by parts} formula holds:
\be
\big(\Dcal(z)\bfm_1\big)\w \bfm_2=\Dcal(z)\big(\bfm_1\w\ovDc(z)\bfm_2\big), \qquad \forall\bfm_1,\bfm_2\in\bw M.\label{eq3:intprtt}
\ee
\qed}
\end{prop}

\claim[{\bf Duality}] Let now $M_0$ be a free abelian group with countable basis $\frakB_0$ as in Section~\ref{preli}.  Let  $\beta_j\in M_0^\vee:=\Hom_\ZZ(M_0,\ZZ)$ such that $\beta_j(b_i)=\delta_{ij}$. The {\em restricted dual} of $M_0$ is
  $M_0^\Vee:=\bigoplus_{j\geq 0}\ZZ\cdot \beta_j$.
The equality
 \[
 \mu_1\w\cdots\w\mu_r(m_1\w\cdots\w m_r):=\det(\mu_i(m_j))_{1\leq i,j\leq r},
 \]
defines a  natural identification of $\bw^rM^\Vee_0$ with $(\bw^rM_0)^\Vee$. In particular
 \[
 (\beta_{\lambda_r}\w\beta_{1+\lambda_{r-1}}\w\cdots\w \beta_{r-1+\lambda_1})_{\blamb\in\Pcal_r}
 \]
  is the basis of $(\bw^rM_0)^\Vee$ dual of $\bw^r\frakB_0$, i.e. $ \beta_{i_1}\w\cdots\w \beta_{i_r}(b_{j_1}\w\cdots\w
 b_{j_r})=\delta_{i_1j_1}\cdots\delta_{i_rj_r}. $

\begin{defin}\label{deftranspo} {\em  The {\em transpose} of $\Dcal(z)\in HS(\wM_0)$ is the module homomorphism
$\Dcal(z)^T:\wM_0^\Vee\sra \wM_0^\Vee[[z]]$, defined by
\[
(\Dcal(z)^T{\bm\eta})(\bfm)={\bm\eta}(\Dcal(z)\bfm), \qquad \forall ({\bm\eta},\bfm)\in \wM_0^\Vee\times \wM_0.
\]
}
\end{defin}

\begin{prop} { \em If  $\Dcal(z)\in HS(\wM_0)$, then $\Dcal(z)^T$ is an $HS$-derivation of $\wM_0^\Vee$.}
\end{prop}
\proof
By definition,  $\Dcal(z)^T\beta(m)=\beta(\Dcal(z)m)$ for all $\beta\in M_0^\star$. As each ${\bm\eta}\in\wM_0^\Vee$ is a sum of
  homogeneous components,  without loss of generality we may assume ${\bm\eta}\in\bw^rM_0^\Vee$, i.e.
$
{\bm\eta}:=\eta_1\w\cdots\w \eta_r
$
for some $\eta_i\in M_0^\Vee$. Thus
\begin{eqnarray*}
\Dcal(z)^T(\eta_1\w\cdots\w \eta_r)(m_1\w\cdots\w m_r)&=&\eta_1\w\cdots\w \eta_r\big(\Dcal(z)(m_1\w\cdots\w m_r)\big)\\
&=&\eta_1\w\cdots\w \eta_r(\Dcal(z)m_1\w\cdots\w\Dcal(z)m_r)\\
&=&\det(\eta_i(\Dcal(z)m_j))=\det(\Dcal(z)\eta_i(m_j))=\\
&=&\Dcal(z)^T\eta_1\w\cdots\w \Dcal(z)^T\eta_r(m_1\w\cdots\w m_r).
\end{eqnarray*}
The unique $HS$-derivation $\widehat{\Dcal}(z)$ on $\wM_0^*$ such that $\widehat{\Dcal}(z)\eta={\Dcal}(z)^T\eta$
coincides with $\Dcal(z)^T$ when evaluated on $\bw^rM_0^\Vee$. Then it must coincide with it and $\Dcal(z)^T\in HS(\wM_0^\Vee)$.\qed

\vspace{-1pt}
\section{Schubert Derivations on $\ZZ{[}X{]}$}\label{secSchuder}

\claim{} \label{secshur1} With the same notation as Section~\ref{preli},  Prop.~\ref{unicitypr} guarantees
the existence of unique  $HS$-derivations
\[
\sigma_+(z):=\sum_{i\geq 0}\sigm_iz^i:\wM_0\sra \wM_0[[z]]\quad \mathrm{and}\quad  \ovsig_+(z)=\sum_{i\geq 0}(-1)^i\ovsig_iz^i:\wM_0\sra \wM_0[[z]]
\]
 such that $\sigm_+(z)b_i=\sum_{j\geq 0}b_{i+j}z^j$ and $\ovsig_+(z)b_i=b_i-b_{i+1}\cdot z$. In particular $\sigm_jb_i=b_{i+j}$ and $\ovsig_jb_i=0$ if $j\geq 2$. They are one the inverse of the other:
 \[
 \sigm_+(z)\ovsig_+(z)=\ovsig_+(z)\sigm_+(z)=1_{\wM_0}.
 \]
We shall call  them {\em Schubert derivations},  in compliance with
the terminology introduced in~\cite{G1,G2}. The motivation comes from
following Pieri-like  formula (\cite[Theorem 2.4]{G1}): \be
\sigma_i\wb^r_\blamb=\sum_\bmu\wb^r_\mu,\qquad (i\geq
0)\label{eq2:pieri} \ee where the sum is taken  over all the
partitions $\bmu\in\Pcal_r$ such that   $\mu_1\geq
\lambda_1\geq\cdots\geq\mu_r\geq \lambda_r$ and $|\bmu|=|\blamb|+i$.
 In addition, a Giambelli-like formula holds (\cite[Formula~(17)]{G1} or~\cite[Theorem 0.1]{LakTh1}): for all $\blamb\in\Pcal_r$
\be
\wb^r_\blamb=\Delta_\blamb(\bsig_+)\wb^r_0:=\det(\sigm_{\lambda_j-j+i})_{1\leq i,j\leq r}\cdot\wb^r_0,\label{eq3:msgmb}
\ee
where by conventions $\sigm_j=0$ if $j<0$.
If $M_{0,n}:=\bigoplus_{j=0}^{n-1}\ZZ b_j$, formula~(\ref{eq2:pieri}) tells us that  $\bw^rM_{0,n}$ is an irreducible representation of
 the cohomology ring $H^*(G_r(\CC^n),\ZZ)$: the latter is in fact  generated as a $\ZZ$-algebra by the special Schubert cycles $c_i(\Qcal_r)$,
  the $i$-th Chern classes of the universal quotient bundle over $G_r(\CC^n)$,  traditionally denoted by $\sigm_i$.  So the
   reason we are using the same notation, more than an abuse, is to emphasize that we are working precisely with the same objects,
    seeing $\bw^rM_0$ as a module over the cohomology (as in~\cite[p.~303]{BottTu}) of the  Grassmannian  $G_r(\CC^\infty)$.
\claim{}\label{secshur2} We similarly define a sort of ``mirror'' of the Schubert derivation $\sigm_+(z)$, namely
\be
\sigm_-(z):=\sum_{i\geq 0}\sigm_{-i}z^{-i}\in HS(\wM_0)
\ee
which, by definition, is the unique $HS$-derivation such that $\sigm_{-j}b_i=b_{i-j}$ if $i\geq j$ and $0$ otherwise. Its inverse in $\End_\ZZ(\wM_0)[[z^{-1}]]$,
\[
\ovsig_-(z)=\sum_{j\geq 0}(-1)^j\ovsig_{-j}z^{-j},
\]
is the unique $HS$-derivation such that $\ovsig_-(z)b_i:=b_{i}-b_{i-1}z^{-1}$ for all $i\geq 0$. In particular, for all $r\geq 0$:
\be
\ovsig_-(z)\wb^r_0=\wb^r_0.
\ee
We note in
 passing that for $j> 0$,  both $\ovsig_{-j}$ and $\sigm_{-j}$ are locally nilpotent, i.e. for all $\bfm\in\wM_0$ there exists $N\in\NN$
  such that $(\sigm_{-j})^N\bfm=0$ (resp. $(\ovsig_{-j})^N\bfm=0$).

\begin{lem} \label{sigpm} {\em For all $r\geq 1$,  let  $(1^r)$ be the partition $(1,\ldots,1)$, with $r$ parts equal to one. Then for
 all $\blamb\in\Pcal_r$ we have
\[
\ovsig_{\pm r}\wb^r_\blamb=\wb^r_{\blamb\pm(1^r)}.
\]
}
\end{lem}
\proof Indeed
 $\ovsig_{\pm r}\wb^r_\blamb$  (Cf. the defining formulas at the beginning of \ref{secshur1} and~\ref{secshur2}) is the coefficient of the monomial  $z^{\pm r}$ in the expansion of
\[
\ovsig_\pm(z)\wb^r_\blamb=\ovsig_\pm(z)(b_{\lambda_r}\w\cdots\w b_{r-1+\lambda_1})
\]
in powers of $z$.
Since  $\ovsig_\pm(z)$ is a $HS$-derivation, we have
\[
\ovsig_\pm(z)(b_{\lambda_r}\w\cdots\w b_{r-1+\lambda_1})=\ovsig_\pm(z)b_{\lambda_r}\w\cdots\w \ovsig_\pm(z)b_{r-1+\lambda_1}
\]
\[
=(b_{\lambda_r}-b_{\lambda_r\pm 1}z^{\pm 1})\w\cdots\w(b_{r-1+\lambda_1}-b_{r-1+\lambda_1\pm 1}z^{\pm 1})
\]
and is so apparent that the coefficient of $z^{\pm r}$ is $\wb^r_{\blamb\pm(1^r)}$, as desired.\qed
\begin{lem} \label{lems+s-} {\em The following equality holds in $(\bw^rM_0)[[z]]$:
\[
\ovsig_+(z)\bfm=(-1)^rz^r\ovsig_-(z)\ovsig_r\bfm.
\]
}
\end{lem}
\proof Without loss of generality, we may assume $\bfm=\wb^r_\blamb$. In this case
\begin{center}
\begin{tabular}{rcll}
$\ovsig_+(z)\wb^r_\blamb$\hskip-3pt & $=$&\hskip-6pt$\ovsig_+(z)(\wbrl)$&(definition of $\wb^r_\blamb$)\\ \\
&$=$&\hskip-3pt$\ovsig_+(z)b_{\lambda_r}\w\cdots\w \ovsig_+(z)b_{r-1+\lambda_1}$&($\ovsig_+(z)\in HS(\bw M_0)$) \\ \\
&$=$&\hskip-3pt$(b_{\lambda_r}-b_{\lambda_r+1}z)\w\cdots\w (b_{r-1+\lambda_1}-b_{r+\lambda_1}z)$& (definition of $\ovsig_-(z)b_j$)\\ \\
&$=$&\hskip-3pt$(-1)^rz^r\left[(b_{\lambda_r+1}-b_{\lambda_r}z^{-1})\w\cdots\w  (b_{r+\lambda_1}-b_{r-1+\lambda_1}z^{-1})\right]$&(highlights $(-1)^rz^r$)\\ \\
&$=$&\hskip-3pt$(-1)^rz^r[\ovsig_-(z)b_{\lambda_r+1}\w\cdots\w\ovsig_-(z)b_{r+\lambda_1}]$&(definition of $\ovsig_-(z)$) \\ \\
&$=$&\hskip-3pt$(-1)^rz^r\ovsig_-(z)(b_{\lambda_r+1}\w\cdots\w b_{r+\lambda_1})$&($\sigm_-(z)\in HS(\wM_0)$)\\  \\
&$=$&\hskip-3pt$(-1)^rz^r\ovsig_-(z)\ovsig_r\wb^r_\blamb$&(\ref{sigpm} applied  to $\ovsig_r$).\,\,\,\,\,\,\,\,\,\,\,\,\,\,\,\,\,\,\,\,\qedd
\end{tabular}
\end{center}
An analogue of Lemma~\ref{lems+s-} holds for $\ovsig_-(z)$ as well, up to an additional hypothesis.
\begin{lem}\label{lems-s+} {\em For all $\wb^r_\blamb\in \bw^rM_0$ such that  $\ell(\blamb)=r$ {\em (}i.e.  $\lambda_r>0${\em )}, then
\be
\ovsig_-(z)\wb^r_\blamb=(-1)^rz^{-r}\ovsig_+(z)\ovsig_{-r}\wb^r_\blamb.\label{eq2:reverse}
\ee
}
\end{lem}

\smallskip
\noindent
Notice that if $\ell(\blamb)<r$, the right hand side of~(\ref{eq2:reverse}) is zero and then~(\ref{eq2:reverse}) fails to be true in general.
\smallskip

\proof Applying the definition of $\ovsig_-(z)$ we eventually arrive to the equality:
\be
\ovsig_-(z)\wb^r_\blamb=\ovsig_-(z)(b_{\lambda_r}\w b_{1+\lambda_{r-1}}\w \cdots\w b_{r-1+\lambda_1})=(b_{\lambda_r}-b_{\lambda_r-1}z^{-1})\w\cdots\w (b_{r-1+\lambda_1}-b_{r-2+\lambda_1}z^{-1}),\label{eq2:inturn}
\ee
by just imitating the first few steps of the proof of Lemma~\ref{lems+s-}. In turn, the left hand side of~(\ref{eq2:inturn})  can be written as
\[
(-1)^rz^{-r}(b_{\lambda_r-1}-b_{\lambda_r}z)\w\cdots\w (b_{r-2+\lambda_1}-b_{r-1+\lambda_1}z)=(-1)^rz^{-r}\ovsig_+(z)b_{\lambda_r-1}\w\cdots\w \ovsig_+(z)b_{\lambda_1+r-2}
\]
i.e., using that $\ovsig_+(z)\in HS(\wM_0)$:
\[
=(-1)^rz^{-r}\ovsig_+(z)(b_{\lambda_r-1}\w\cdots\w b_{r-1+\lambda_1})=(-1)^rz^{-r}\ovsig_+(z)\ovsig_{-r}\wb^r_\blamb,
\]
where the last equality holds because of the hypothesis $\lambda_r>0$.\qed

Recall Definition~\ref{deftranspo}. The next easy  lemma identifies the transpose $\sigm_+(z)^T$ of $\sigm_+(z)$, relating it with the Schubert derivations on $\wM_0^\Vee$.

\begin{lem}\label{trspbeta} {\em Let  $\sum_{j\geq 0}\sigm_j^Tz^{-j}:=\sigm_-(z)^T\in \End_A(\wM_0^\Vee)[[z^{-1}]]$.
Then $\sigm_{-j}^T\beta_i=\beta_{i+j}$.
}
\end{lem}
\proof In fact, for all $k\geq 0$:
\[
\sigm_{-j}^T\beta_i(b_k)=\beta_i(\sigm_{-j}b_k)=\beta_i(b_{k-j})=\delta_{i,k-j}=\delta_{i+j,k}=\beta_{i+j}(b_k)
\]
which proves the claim.\qed

\vspace{-1pt}
\section{Proof of Theorem~\ref{thm01}}\label{EGC}

The starting point is the following well--known criterion (\cite[Section 4]{Arb} or, in an infinite dimensional context,  \cite[Proposition 7.2]{KR}).

\begin{prop} {\em An element $\bfm\in\bw^rM_0$ belongs to
$\Gcal_r$ if and only if the equality \be \sum_{i\geq 0} (\beta_i\lrcorner \bfm)\otimes(b_i\w
\bfm)=0\label{eq4:critkac} \ee holds
in $\bw^{r-1}M_0\otimes\bw^{r+1}M_0$.} \end{prop}
\proof See e.g.~\cite[Theorem~6.1.7]{GatSal2}.

\smallskip
Recall from Section~\ref{secshur1} and Lemma~\ref{trspbeta}
  that
  \[
  \sum_{i\geq 0}b_iz^i=\sigm_+(z)b_0\qquad \mathrm{and}\qquad \sum_{j\geq 0}\beta_jz^{-j-1}=z^{-1}\sigm_-(z)^T\beta_0.
  \]
 Equation~(\ref{eq4:critkac}) can be rewritten, imitating~\cite{Arb,KR},
  in the   equivalent form
\be
\Res_z(z^{-1}\sigm_-(z)^T\beta_0)\lrcorner\bfm\otimes \sigm_+(z)b_0\w\bfm ) )=0,\label{eq:decres}
\ee
i.e.  $\bfm\in \Gcal_r$ if and only it satisfies equation~(\ref{eq:decres}). We have:

\begin{prop} { \em The following equality holds in
$\bw^{r+1}M_0${\em :} \be \sigm_+(z)b_0\w\bfm
=(-1)^rz^r\sigm_+(z)\ovsig_-(z)(b_0\w
\ovsig_r\bfm).\label{eq4:bzweta} \ee }
 \end{prop}

\proof First of all 
\be 
\sigm_+(z)b_0\w\bfm
=\sigm_+(z)\big(b_0\w
\ovsig_+(z)\bfm\big),\label{eq4:bzwetaa} 
\ee
 because of integration by parts~(\ref{eq3:intprtt}). Lemma~\ref{lems+s-} applied to $\ovsig_+(z)\bfm$ gives, after
simplification:
\[
b_0\w \ovsig_+(z)\bfm=(-1)^rz^r\ovsig_-(z)(b_0\w\ovsig_r\bfm).
\]
Substituting in the last side of~(\ref{eq4:bzwetaa}) gives~(\ref{eq4:bzweta}), as desired.\qed
\begin{prop}  { \em The following equality holds in $\bw^{r-1}M_0$:
\be
(z^{-1}\sigm_-(z)^T\beta_0)\lrcorner \bfm =(-1)^{r-1}z^{-r}\ovsig_+(z)\ovsig_{-r+1}(\beta_0\lrcorner \sigm_-(z)\bfm).\label{eq4:bzwetax}
\ee
}
\end{prop} 
\proof
Let $\gamma\in\bw^{r-1}M_0^\Vee$ be arbitrarily chosen. Then
\begin{center}
\begin{tabular}{rrll}
$\gamma(z^{-1}\sigm_-(z)^T\beta_0\lrcorner \bfm)$&$=$&$z^{-1}(\sigm_-(z)^T\beta_0\w\gamma)\bfm$&\hskip 20pt  (by definition~(\ref{eq4:contra}))\\ \\
&$=$&$z^{-1}\sigm_-(z)^T(\beta_0\w \ovsig_-(z)^T\gamma)\bfm$&\hskip 20pt (integration by parts~(\ref{eq3:intprtt}))\\ \\
&$=$&$z^{-1}(\beta_0\w \ovsig_-(z)^T\gamma)\sigm_-(z)\bfm$&\hskip 20pt (definition of $\sigm_-(z)^T$)\\ \\
&$=$&$z^{-1}\ovsig_-(z)^T\gamma(\beta_0\lrcorner \sigm_-(z)\bfm)$&\hskip 20pt (definition~(\ref{eq4:contra}) of contraction) \\ \\
&$=$&$z^{-1}\gamma(\ovsig_-(z)^T(\beta_0\lrcorner \sigm_-(z)\bfm)$&\hskip20pt (definition of $\ovsig_-(z)^T$)
\end{tabular}
\end{center}
whence the equality
\[
z^{-1}\sigm_-(z)^T\beta_0\lrcorner \bfm=z^{-1}\ovsig_-(z)(\beta_0\lrcorner \sigm_-(z)\bfm).
\]
Notice now that $\beta_0\lrcorner \ovsig_-(z)\bfm$ is a linear combination of elements $\wb^{r-1}_\blamb$ associated to partitions
 of length exactly $r-1$. Thus we can apply Lemma~\ref{lems-s+} to get
\[
\hskip62pt z^{-1}\sigm_-(z)^T\beta_0\lrcorner \bfm=z^{-1}\ovsig_-(z)(\beta_0\lrcorner \sigm_-(z)\bfm)=(-1)^{r-1}z^{-r}\ovsig_+(z)\ovsig_{-r+1}(\beta_0\lrcorner \sigm_-(z)\bfm).\hskip28pt \qed
\]
\claim{} Substitution of  expressions ~(\ref{eq4:bzwetax})  and~(\ref{eq4:bzwetaa}) into~(\ref{eq:decres})  concludes the proof of  Theorem~\ref{thm01}.\qed

\vspace{-1pt}
\section{The Grassmann Cone in a Polynomial Ring}\label{GCPR}

\claim{} \label{isophir} Notation and convention as in Section~\ref{preli}. 
Let $\bw^rM_r$ be the  $B_r$-module structure on $\bw^rM_0$ given by:
\[
e_i\wb^r_\blamb:=\ovsig_i\wb^r_\blamb,
\]
which turns   $e_i\in B_r$  into an eigenvalue of $\ovsig_i$. Accordingly,  we have
\[
\ovsig_+(z)\wb^r_\blamb=E_r(z)\wb^r_\blamb
\]
and then
\begin{eqnarray*}
\sigm_+(z)\wb^r_\blamb&=&\sigm_+(z)\left(E_r(z){\sum_{n\geq 0}h_nz^n}\right)\wb^r_\blamb=\\
&=&\sigm_+(z)\ovsig_+(z)(\sum_{n\geq
0}h_nz^n)\wb^r_\blamb=\sum_{n\geq 0}h_n\wb^r_\blamb z^n.
\end{eqnarray*}
Thus  $h_i$ is the eigenvalue of $\sigm_i$ seen as endomorphism of $\bw^rM_r$, i.e.
$
\sigm_i\wb^r_\blamb=h_i\wb^r_\blamb.
$
In particular, using~(\ref{eq3:msgmb}),  the homomorphism of abelian groups  $\phi_r:B_r\sra \bw^rM_0$ given by
\be
\Delta_\blamb(H_r)\mapsto \Delta_\blamb(H_r)\wb^r_0=\Delta_\blamb(\bsig_+(z))\wb^r_0=\wb^r_\blamb,\label{eq8:bfcor}
\ee
 is an isomorphism, as it maps the $\ZZ$-basis $(\Delta_\blamb(H_r))$  of $B_r$ to the  $\ZZ$-basis $\wb^r_\blamb$ of $\bw^rM_0$. By abuse of notation, for all $\bfm\in\bw^rM_0$ we shall  write
 \[
{\bfm\over \wb^r_0}:=\phi_r^{-1}(\bfm),
 \]
 i.e. for the unique element of $B_r$ carrying the $B_r$-basis element $\wb^r_0$ of $\bw^rM_0$ to $\bfm$.

\begin{defin} {\em Let $\sigm_-(z),\ovsig_-(z):B_r\sra B_r[z^{-1}]$ be defined as:
\[
\sigm_-(z)\Delta_\blamb(H_r):={\sigm_-(z)\wb^r_\blamb\over \wb^r_0}\qquad \mathrm{and}\qquad \ovsig_-(z)\Delta_\blamb(H_r):={\ovsig_-(z)\wb^r_\blamb\over \wb^r_0}.
\]
}
\end{defin}
The $\sigm_-(z)$-image of $h_n=\Delta_{(n)}(H_r)$  could in principle depend on the integer $r$. However this is not the case.

\begin{prop}\label{prophnz} {\em For all $r\geq 1$, the following equalities hold in the ring $B_r[z^{-1}]${\em :}
\be
\sigm_-(z)h_n=\sum_{j=0}^nh_{n-i}z^{-j}\qquad \mathrm{and}\qquad \ovsig_-(z)h_n=h_{n}-h_{n-1}z^{-1}. \label{eq5:s.zhn}
\ee
}
\end{prop}

\proof
We have:
\smallskip
\begin{center}
\begin{tabular}{rlll}
$(\sigm_-(z)h_n)\wb^r_0$&$=$&$\sigm_-(z)(h_n\wb^r_0)$&(definition of $\sigm_-(z)h_n$)\\ \\
&$=$&$\sigm_-(z)(\wb^{r-1}_0\w b_{r-1+n}))$&(writing $\wb^r_{(n)}$ as $\wb^{r-1}_0\w b_{r-1+n}$) \\  \\
&$=$&$\sigm_-(z)\wb^{r-1}_0\w \sigm_-(z) b_{r-1+n}$&(since $\sigm_-(z)\in HS(\wM_0)$)\\ \\
&$=$&$\sum_{j=0}^{r-1+n}\left(\wb^{r-1}_0\w b_{r-1+n-j}z^{-j}\right)$&(apply~(\ref{eq2:HSDA}) and the definition of $\sigm_-(z)$)\\ \\
&$=$&$(\sum_{j=0}^{n}h_{n-j}z^{-j})\wb^r_0$,&
\end{tabular}
\end{center}
which proves  the first equality in~(\ref{eq5:s.zhn}).
To prove the second equality of~(\ref{eq5:s.zhn}),  we can argue either by observing that $\ovsig_-(z)$ is the inverse of $\sigm_-(z)$ or again by direct computation:
\begin{eqnarray*}
\hskip37pt (\ovsig_-(z)h_n)\wb^r_0&=&\ovsig_-(z)(h_n\wb^r_0)=\ovsig_-(z)\wb^r_n=\ovsig_-(z)\wb^{r-1}_0\w \ovsig_-(z)b_{r-1+n}\\
&=&\wb^{r-1}_0\w (b_{r-1+n}-b_{r-1+n-1}z^{-1})=(h_n-h_{n-1}z^{-1})\wb^r_0.\hskip36pt \qed
\end{eqnarray*}

\begin{prop}\label{prop1stvo} {\em We have:
\be
b_0\w \ovsig_r\wb^r_\blamb=\wb^{r+1}_\blamb=\Delta_\blamb(H_{r+1})\wb^{r+1}_0.\label{eq:eqsab}
\ee
}
\end{prop}

\proof
Indeed
\begin{center}
\begin{tabular}{rlll}
$b_0\w \ovsig_r\wb^r_\blamb$&$=$&$b_0\w \ovsig_r(b_{\lambda_{r}}\w b_{1+\lambda_{r-1}}\w\cdots\w b_{1+\lambda_r})$&(definition of $\wb^r_\blamb$)\\ \\
&$=$&$b_0\w b_{1+\lambda_{r}} \w b_{2+\lambda_{r-1}}\w\cdots\w
b_{r+\lambda_r}$&(Lemma~\ref{sigpm})\\ \\
&$=$&$\wb^{r+1}_\blamb$&(definition of $\wb^{r+1}_\blamb$)
\end{tabular}
\end{center}
In the $B_{r+1}$-module $\bw^{r+1}M_r$ we have then equality~(\ref{eq:eqsab}) (due to  $\phi_{r+1}(\Delta_\blamb(H_{r+1}))=\wb^{r+1}_\blamb$, by~(\ref{eq8:bfcor})).\qed


In other words the expressions of $\sigm_-(z)h_n$ and $\ovsig_-(z)h_n$ in the ring $B_r[z^{-1}]$ do not depend on the integer $r$.

\claim{} Let us agree to simply write $\Res(f)$ for the coefficient of $X^{-1}$ of a formal Laurent series $f\in\ZZ((X^{-1}))$. To prove Theorem~\ref{0stmnthm} below, we need a powerful result due to Laksov and Thorup~\cite[Theorem 0.1.(2)]{LakTh1} (see also~\cite{LakTh2}).  Let
  us introduce a few new pieces of notation.
Let
\[
p_r(X)=X^rE_r\left({1\over X}\right)=X^r-e_1X^{r-1}+\cdots+(-1)^re_r\in B_r[X],
\]
be the generic polynomial of degree $r$.
If $f:=a_0X^\lambda+a_1X^{\lambda-1}+\cdots +a_\lambda\in B_r[X]$ is any polynomial of degree $\leq \lambda$, then   an easy computation shows
that
\be
\Res\left({X^{i-1}f(X)\over p_r(X)}\right):=\Res\,{X^{i-1}f(X)\over X^r}\left(1+{h_1\over X}+{h_2\over X^2}+\cdots\right)=\sum_{j=0}^\lambda a_jh_{i-r+\lambda-j}.\label{eq7:easy}
\ee
for all $i\geq 0$.
Following~\cite{LakTh1}, the {\em residue} of
 $(f_0,f_1,\ldots,f_{r-1})\in B_r[X]^r$ is,  by definition:
\be
\Res(f_0,f_1,\ldots,f_{r-1})=\left|\begin{matrix}\Res(f_0)&\Res(f_1)&\cdots&\Res(f_{r-1})\cr
\Res(Xf_0)&\Res(Xf_1)&\cdots&\Res(Xf_{r-1})\cr
\vdots&\vdots&\ddots&\vdots\cr
\Res(X^{r-1}f_0)&\Res(X^{r-1}f_1)&\cdots&\Res(X^{r-1}f_{r-1})
\end{matrix}\right|.\label{eq:reslkth}
\ee
We shall use the following
\begin{thm}[\cite{LakTh1}, Theorem 0.1 (2)] \label{thmlakth1} {\em Let $f_0,f_1,\ldots, f_{r-1}\in B_r[X]$. Then
\begin{eqnarray}
&& f_0(\sigm_1)b_0\w f_1(\sigm_1)b_0\w\cdots\w
f_{r-1}(\sigm_1)b_0\cr\cr &=& \Res\left({f_{r-1}(X)\over
p_r(X)}, {f_{r-2}(X)\over p_r(X)}, \ldots, {f_{0}(X)\over
p_r(X)}\right) b_0\w b_1\w\cdots\w b_{r-1}.\hskip25pt\cr\cr
&=&\det\left(\Res \left({X^{i-1}f_{r-j}(X)\over p_r(X)}\right)\right)_{1\leq i,j\leq r}b_0\w b_1\w\cdots\w b_{r-1}.
\label{eq:detflakth}
\end{eqnarray}
}
\end{thm}
Denote by  $\sigm_-(z)H_r$ and $\ovsig_-(z)H_r$, respectively, 
 the following $B_r[z^{-1}]$-valued sequences:
 \[
 (\sigm_-(z)h_j)_{j\in\ZZ}=(\sum_{i=0}^jh_{j-i}z^{-i})_{j\in\ZZ}\qquad{\mathrm and}\qquad (\ovsig_-(z)h_j)_{j\in\ZZ}=(h_j-h_{j-1}z^{-1})_{j\in\ZZ}.
 \]

\begin{thm} \label{0stmnthm}{\em The operators $\sigm_-(z),\ovsig_-(z):B_r\sra B_r[z^{-1}]$ commute with taking  Schur determinants, i.e.{\em :}
\be
\ovsig_-(z)\Delta_\blamb(H_r)=\Delta_\blamb(\ovsig_-(z)H_r)\qquad\mathrm{and}\qquad \sigm_-(z)\Delta_\blamb(H_r)=\Delta_\blamb(\sigm_-(z)H_r).\label{eq:commsd}
\ee
}
\end{thm}

\proof To prove the first of equalities~(\ref{eq:commsd}), one observes that
\begin{eqnarray*}
(\ovsig_-(z)\Delta_\blamb(H_r))\wb^r_0&=&\ovsig_-(z)\wb^r_\blamb=\ovsig_-(z)b_{\lambda_r}\w\cdots\w\ovsig_-(z)b_{r-1+\lambda_1}\\
&=&f_0(\ovsig_{1})b_0\w f_1(\sigm_1)b_0\w\cdots\w
f_{r-1}(\sigm_1)b_0,
\end{eqnarray*}
where $f_i(\ovsig_{1})$ stands  for $f_i(X)=X^{\lambda_{r-i}}-X^{\lambda_{r-i}-1}z^{-1}$ evaluated at $X=\sigm_{1}$.
By Theorem~\ref{thmlakth1} and formula~(\ref{eq7:easy}):
\begin{eqnarray*}
f_0(\sigm_{1})b_0\w f_1(\sigm_{1})b_0\w\cdots\w f_{r-1}(\sigm_{1})b_0&=&\Res\left({f_{r-1}(X)\over p_r(X)},\ldots,{f_{0}(X)\over p_r(X)}\right)b_0\w b_1\w\cdots\w b_{r-1}\\
&=&\det(h_{\lambda_j-j+i}-h_{\lambda_j-j+i-1}z^{-1})_{1\leq i,j\leq r}b_0\w b_1\w\cdots\w b_{r-1}\\
&=&\Delta_\blamb(\sigm_-(z)H_r)b_0\w b_1\w\cdots\w b_{r-1}.
\end{eqnarray*}
To verify the second equality of~(\ref{eq:commsd}) we exploit the first one, just proven, and the fact that $\sigm_-(z)$ and $\ovsig_-(z)$ are one inverse of the other. Then
\[
\hskip19pt \sigm_-(z)\Delta_\blamb(H_r)=\sigm_-(z)\Delta_\blamb(\ovsig_-(z)\sigm_-(z)H_r)=\sigm_-(z)\ovsig_-(z)\Delta_\blamb(\sigm_-(z)H_r)=\Delta_\blamb(\sigm_-(z)H_r).\hskip 19pt\qed
\]

\begin{lem} {\em  For all $\blamb\in\Pcal_r$, the following
equality holds in $\bw^{r-1}M_0$ \be \ovsig_{-r+1}(\beta_0\lrcorner
\wb^r_\blamb)=\Delta_\blamb(H_{r-1})\wb^{r-1}_0.\label{eq8:lemprevvv}
\ee
}
\end{lem}

\proof If $\ell(\blamb)=r$ then $\lambda_r>0$ and then $\beta_0\lrcorner\wb^r_\blamb=0$. So, both sides of~(\ref{eq8:lemprevvv})
 vanish.  If $\ell(\blamb)<r$ instead, one may write
\[
\wb^r_\blamb=b_0\w \wb^{r-1}_{\blamb+(1^{r-1})},
\]
where $\blamb+(1^{r-1})=(\lambda_1+1,\ldots,\lambda_{r-1}+1)$. Thus
\[
\ovsig_{-r+1}(\beta_0\lrcorner (b_0\w \wb^{r-1}_{\blamb+(1^{r-1})}))=\ovsig_{-r+1}(\wb^{r-1}_{\blamb+(1^{r-1})})=\wb^{r-1}_\blamb=\Delta_\blamb(H_{r-1})\wb^{r-1}_0
\]
 as claimed. \qed

\begin{prop} \label{prop2stvo}  {\em The equality below holds in $\bw^{r-1}M_0${\em :}
\[
\ovsig_{-r+1}(\beta_0\lrcorner \sigm_-(z)\wb^r_\blamb)=\Delta_\blamb(\sigm_-(z)H_{r-1})\wb^{r-1}_0
\]
}
\end{prop}
\proof
If one defines  $\Delta_{(n)}(\sigm_-(z)H_0)$ to be $z^{-n}$, the equality holds for $r=1$. To  check the formula in the remaining cases, 
let us preliminary observe that
$
\Delta_\blamb(\sigm_-(z)H_r)
$
is a linear combination $\sum_{\bmu\subseteq \blamb}a_\bmu(z^{-1})\Delta_\bmu(H_r)$, whose coefficients  $a_\bmu(z^{-1})\in \ZZ[z^{-1}]$
   do not depend on the chosen $r>1$. Thus:
 
\pagebreak
\begin{eqnarray*}
\hskip32pt \ovsig_{-r+1}\left(\beta_0\lrcorner \sigm_-(z)\wb^r_\blamb\right)&=&\ovsig_{-r+1} \big(\beta_0\lrcorner (\sigm_-(z)(\Delta_\blamb(H_r)\wb^r_0)\big)\cr\cr
&=&\ovsig_{-r+1}\big(\beta_0\lrcorner (\sigm_-(z)\Delta_\blamb(H_r))\wb^r_0)\big)\\ \\
&=&\ovsig_{-r+1}\big(\beta_0\lrcorner (\Delta_\blamb(\sigm_-(z) H_r)\wb^r_0\big) \\ \\
&=&\ovsig_{-r+1}\beta_0\lrcorner\sum_{\bmu\subseteq \blamb}a_\bmu(z^{-1})\Delta_\bmu(H_r)\wb^r_0\cr
&=&\sum_{\bmu\subseteq \blamb}a_\bmu(z^{-1})\cdot \ovsig_{-r+1}(\beta_0\lrcorner\Delta_\bmu(H_r)\wb^r_0)\cr\cr
&=&\sum_{\bmu\subseteq \blamb}a_\bmu(z^{-1})\cdot \Delta_\blamb(H_{r-1})\wb^{r-1}_0=\Delta_\blamb(\sigm_-(z)H_{r-1})\wb^{r-1}_0.\hskip30pt \qed
\end{eqnarray*}

\claim{} Let $\Gamma_r(z):B_r\sra B_{r+1}((z))$ and $\Gamma^\Vee_r(z):B_r\sra B_{r-1}((z))$ defined by:
\be
\Gamma_r(z)\Delta_\blamb(H_r):={1\over E_{r+1}(z)}\ovsig_-(z)\Delta_\blamb(H_{r+1}),\label{eq8:gmv}
\ee
and
\be
 \Gamma^\Vee_r(z)\Delta_\blamb(H_r):=E_{r-1}(z)\Delta_\blamb(\sigm_-(z)H_{r-1}).\label{eq8:gmvee}
\ee
Expression~(\ref{eq8:gmv}) can be also written in the form
\[
\Gamma_r(z)\Delta_\blamb(H_r)={1\over E_{r+1}(z)}\Delta_\blamb(\ovsig_-(z)H_{r+1}),
\]
by virtue of Theorem~\ref{0stmnthm}, according which   $\ovsig_-(z)\Delta_\blamb(H_{r+1})=\Delta_\blamb(\ovsig_-(z)H_{r+1})$ for all
 $\blamb\in\Pcal_{r+1}$. Similarly, the equality
$
\Delta_\blamb(\sigm_-(z)H_{r-1})= \sigm_-(z)\Delta_\blamb(H_{r-1})
$
surely  holds for all  $\blamb\in\Pcal_{r-1}$. However, if  $\ell(\blamb)=r$, in general $\Delta_\blamb(\sigm_-(z)H_{r-1})\neq  \sigm_-(z)\Delta_\blamb(H_{r-1})$, because $\Delta_\blamb(\sigm_-(z)H_{r-1})\neq 0$ in spite of the vanishing of  $\Delta_\blamb(H_{r-1})$.  For  example,  if
  $r=1$,  $\sum_{n\geq 0}h_nz^n=(1-e_1z)^{-1}=\sum_{n\geq 0}e_1^nz^n$, i.e. $h_n=h_1^n$. This  implies
\[
\Delta_{(1,1)}(H_1)=\left|\begin{matrix}h_1&1\cr h_2&h_1\end{matrix}\right|=h_1^2-h_2=h_1^2-h_1^2=0.
\]
On the other hand
\[
\Delta_{(1,1)}(\sigm_-(z)H_1))=\left|\begin{matrix}h_1+\displaystyle{1\over z}&1\cr\cr h_2+\displaystyle{h_1\over z}+\displaystyle{1\over z^{2}}&h_1+\displaystyle{1\over z}\end{matrix}\right|=\left(h_1+{1\over z}\right)^2-h_2-{h_1\over z}-{1\over z^2}={h_1\over z}\neq 0.
\]


\noindent
{\bf Proof of Theorem~\ref{thm02}.} According to Theorem~\ref{thm01},  it follows that $\sum_{\blamb\in\Pcal_{r,n}}a_\blamb\wb^r_\blamb\in \Gcal_r$ if and only if
\be
\Res_z \sum_{\blamb,\bmu\in\Pcal_r}a_\blamb a_\bmu \ovsig_+(z)\ovsig_{-r+1}(\beta_0\lrcorner \sigm_-(z)\wb^{r-1}_\blamb)\otimes_\ZZ\sigm_+(z)\ovsig_-(z)(b_0\w\ovsig_r\wb^{r+1}_\bmu)\label{eq7:eigenv}
\ee
vanishes  in $(\bw^rM_{r-1}\otimes\bw^rM_{r+1})((z))$.
Since $\ovsig_{-r+1}(\beta_0\lrcorner\sigm_-(z)\wb^{r-1}_0)\in \bw^{r-1}M_{r-1}$, it is an eigenvector of
 $\ovsig_+(z)$ corresponding to the eigenvalue $E_{r-1}(z)$. Similarly, $\ovsig_-(z)(b_0\w \ovsig_{-r}\wb^r_\blamb)$ belongs to $\bw^{r+1}M_0$, 
  which is an eigenmodule of $\sigm_+(z)$ corresponding to the  eigenvalue $1/E_{r+1}(z)$. Thus,   expression~(\ref{eq7:eigenv}) can be rewritten as
\begin{eqnarray*}
0&=&\Res_z\, \sum_{\blamb,\bmu\in\Pcal_{r}}a_\blamb a_\bmu E_{r-1}(z)(\ovsig_{-r+1}(\beta_0\lrcorner \sigm_-(z)\wb^r_\blamb)\otimes {1\over E_{r+1}(z)}\ovsig_-(z)(b_0\w\ovsig_r\wb^{r+1}_\blamb)=\\
&=&\Res_z\sum_{\blamb,\bmu\in\Pcal_{r}} E_{r-1}(z)\Delta_\blamb(\sigm_-(z)H_{r-1})\wb^{r-1}_0\otimes {1\over E_{r+1}(z)}\ovsig_-(z)\Delta_\bmu(H_{r+1})\wb^{r+1}_0\\
&=&\left(\Res_z\sum_{\blamb,\bmu\in\Pcal_{r}}
E_{r-1}(z)\Delta_\blamb(\sigm_-(z)H_{r-1})\wb^{r-1}_0\otimes {1\over
E_{r+1}(z)}\ovsig_-(z)\Delta_\bmu(H_{r+1})\right)\wb^{r-1}_0\otimes
\wb^{r+1}_0,
\end{eqnarray*}
where in last equality we used Propositions~\ref{prop1stvo} and~\ref{prop2stvo}.
Let us now identify the tensor product  $B_{r-1}\otimes_\ZZ B_{r+1}$ with the polynomial ring in $2r$ indeterminates
\be
B_{r-1}\otimes B_{r+1}=\ZZ[e_1' ,\ldots, e'_{r-1},e_1'',\ldots,e_{r+1}''].\label{eq:rngid}
\ee
Then~(\ref{eq1:mnthm}) follows because  $\wb^{r-1}_0\otimes \wb^{r+1}_0$ is a basis of $\bw^{r-1}M_{r-1}\otimes_\ZZ\bw^{r+1}M_{r+1}$ over $B_{r-1}\otimes_\ZZ B_{r+1}$ though of a ring through the indentification~(\ref{eq:rngid}).\qed

\claim{} Let $E'_{r-1}(z)=1-e'_1z+\cdots+(-1)^{r-1}e'_{r-1}z^{r}\in
B'_{r-1}[z]$ and
$E''_{r+1}(z)=1-e''_1z+\cdots+(-1)^{r+1}e''_{r+1}z^{r+1}\in
B'_{r+1}[z]$. Similarly, let
\[
H'_{r-1}(z)=\sum_{n\geq 0}h'_nz^n\qquad \mathrm{and}\qquad
H''_{r-1}(z)=\sum_{n\geq 0}h''_nz^n
\]
 be the inverses of $E'_{r-1}(z)$ and $E''_{r+1}(z)$ in $B'_{r-1}[[z]]\cong B_{r+1}[[z]]$ and $B''_{r+1}[[z]]$
respectively. Then formula~(\ref{eq1:mnthm})  reads
\be
\Res_z {E'_{r-1}(z)\over E''_{r+1}(z)}\ \sum_{\blamb,\bmu}a_\blamb
a_\bmu\Delta_\blamb(\sigm_-(z)H'_{r-1})\cdot
\ovsig_-(z){\Delta_\bmu(H''_{r+1})}=0.\label{eq:inthtp}
\ee

\vspace{-1pt}
\section{An Example} \label{SecEx}
To show how formula~(\ref{eq:inthtp}) works, in this section we proceed  to determine all the polynomials
\begin{eqnarray}
p(H_2)=a_0+a_1h_1+a_2h_2+ a_{11}\Delta_{(11)}(H_2)+a_{21}\Delta_{(21)}(H_2)+a_{22}\Delta_{(22)}(H_2)\in B_2,\label{eq:klqdvo}
\end{eqnarray}
such that $p(H_2)\wb^2_0$ belongs to the Grassmann cone $\Gcal_2$ of $\bw^2M_0$. Needless to say, we expect to find the same
 expression of the Klein quadric hypersurface of $\PP^5$. For this purpose we need a few preliminaries to speed up
  computations, which by the way could be carried out directly.

\claim{} Consider the two sub-modules  $\sigm_nM_0:=\bigoplus_{i\geq 0}\ZZ\cdot b_{n+i}$ and $M_{0,n}:=\bigoplus_{i=0}^{n-1} \ZZ b_i$ of $\bw^rM_0$.
 They fits in the exact sequence
\[
0\lra \sigm_nM_0\lra M_0\stackrel{t_n}{\lra} M_{0,n}\lra 0,
\]
where $t_n$ is the truncation which maps $b_i$ to itself if $0\leq i\leq n-1$ and to zero otherwise.
 Hence  $\bw^rM_{0,n}$ is the submodule $\bigoplus_{\blamb\in\Pcal_{r,n}}\ZZ\wb^r_\blamb$ of $\bw^rM_0$. It can be seen itself as
  the epimorphic image of the truncation   $t^r_n:\bw^rM_0\sra \bw^rM_{0,n}$ which maps  $\wb^r_\blamb$ to itself if
   $\blamb\in\Pcal_{r,n}$ and to $0$ otherwise. Again, we have an exact sequence
 \[
0\lra \sigm_nM_0\w\bw^{r-1}M_0\lra \bw^rM_0\lra \bw^rM_{0,n}\lra 0,
\]
 where  $\sigm_nM_0\w\bw^{r-1}M_0$ is the submodule of $\bw^rM_0$ generated by all $b_{i_1}\w\cdots\w b_{i_r}\in \bw^{r}M_0$ such
  that $i_j\geq n$ for at least one $1\leq j\leq r$.

Let $I_{r,n}$ be the kernel of the composition  $t^r_n\circ \phi_r:B_r\sra \bw^rM_{0,n}$, where  $\phi_r:B_r\sra \bw^rM_0$ is the
 isomorphism~(\ref{eq8:bfcor}). It is a simple exercise to check that  $h_{n-r+1+j}\in I_{r,n}$ for all $j\geq 0$. Indeed,
   a simple argument (see~\cite[Chapter 5]{GatSal2} or~\cite{G2} for details)
shows that $I_{r,n}=(h_{n+r-1},\ldots, h_n)$. Let
\[
\pi_n:B_r\lra B_{r,n}:={\ZZ[e_1,\ldots, e_r]\over (h_{n-r+1},\ldots, h_n)}
\]
be the canonical projection. Then
$
B_{r,n}=\bigoplus_{\blamb\in\Pcal_{r,n}}\ZZ\Delta_\blamb(H_{r,n})
$,
where by $H_{r,n}$ we have denoted  the sequence  $(h_j+I_{r,n})_{j\geq 0}$  of $B_{r,n}$ and then
 $\Delta_\blamb(H_{r,n})=\Delta_\blamb(H_r)+I_{r,n}$.
Clearly, $\Delta_\blamb(H_{r,n})\wb^r_0=\Delta_\blamb(H_{r})\wb^r_0$ if $\blamb\in\Pcal_{r,n}$. Define
 $\sigm_{-i}(h_j+I_{r,n})=\sigm_jh_n+I_{r,n}$ (resp.
 $\ovsig_{-i}(h_j+I_{r,n})=\ovsig_jh_n+I_{r,n}$).
Then Theorem~\ref{thm02} has the following

\begin{corol}\label{corprtr} {\em  A tensor $\sum_{\blamb\in\Pcal_{r,n}}a_\blamb\wb^r_\blamb$ is decomposable if and only if
\[
\Res_{z}\sum_{\blamb,\bmu\in\Pcal_{r,n}}a_\blamb a_\bmu E_{r-1}(z)\Delta_\blamb(\sigm_-(z)H_{r-1,n})\otimes {1\over E_{r+1}(z)}\Delta_\blamb(\ovsig_-(z)H_{r+1,n})=0
\]
in the tensor product $(B_{r-1,n}\otimes_\ZZ B_{r+1,n})((z))$.
}
\end{corol}

\proof
In fact by Theorem~\ref{thm02},  $\sum_{\blamb\in\Pcal_{r,n}}a_\blamb\wb^r_\blamb\in \Gcal_r$ if and only if
\begin{eqnarray*}
\hskip19pt 0&=&\Res_z\sum_{\blamb,\bmu\in\Pcal_{r}} E_{r-1}(z)\Delta_\blamb(\sigm_-(z)H_{r-1})\wb^{r-1}_0\otimes {1\over E_{r+1}(z)}\ovsig_-(z)\Delta_\bmu(H_{r+1})\wb^{r+1}_0\\ 
&=&\Res_z\sum_{\blamb,\bmu\in\Pcal_{r}}
E_{r-1}(z)\Delta_\blamb(\sigm_-(z)H_{r-1,n})\wb^{r-1}_0\otimes
{1\over E_{r+1}(z)}\ovsig_-(z)\Delta_\bmu(H_{r+1,n})\wb^{r+1}_0.\hskip20pt  \qed
\end{eqnarray*}

Moreover formula~(\ref{eq:inthtp}) and  Corollary~\ref{corprtr} easily imply:
\begin{corol} \label{cortrunc} {\em The  tensor
$
\sum_{\blamb\in \Pcal_{r,n}}a_\blamb\wb^r_\blamb\in\bw^rM_0
$
is decomposable  if and only if
\[
\hskip66pt \Res_z \,{E'_{r-1}(z)\over E''_{r+1}(z)} \hskip-2pt
\sum_{\blamb,\bmu\in \Pcal_{r,n}}\hskip-7pt a_\blamb
a_\bmu\Delta_\blamb(\sigm_-(z)H'_{r-1,n})\cdot
\ovsig_-(z){\Delta_\bmu(H''_{r+1,n})}=0.\hskip47pt \qed
\]
} \end{corol}

\vspace{-12pt}

\claim{}\label{explug24} Going back to our original purpose of finding all $p(H_2)\in B_2$ such that $p(H_2)\wb_0^2\in\Gcal_2$, where
 $p(H_2)$ is like in~(\ref{eq:klqdvo}),
  we use the
$B_2$-module structure of $\bw^2\hskip-2ptM_{2}:=B_2\otimes
\bw^2\hskip-2ptM_0$. Since $p(H_2)\wb^2_0\in
\bw^2\hskip-2ptM_{2,4}$,  the $B_2$-module structure of
$\bw^2\hskip-2ptM_2$  factorizes  through that of $B_{2,4}$, i.e. $p(H_2)\wb^2_0=p(H_{2,4})\wb^2_0$.
We have
\[
B_{1,4}:={B_1\over (h_4)}\cong {\ZZ[x]\over (x^4)},\qquad \mathrm{and}\qquad B_{3,4}:={B_3\over (h_2,h_3,h_4)}\cong {\ZZ[y]\over (y^4)},
\]
where we have set $x=e_1+(h_4)$ and $y=e_1+(h_2,h_3,h_4)$. In particular,  $h_i+(h_4)=x^i$ and $h_1+(h_2,h_3,h_4)=y$.
Thus
\begin{eqnarray*}
p(\sigm_-(z)H_{1,4})&=&a_0+a_1\sigm_-(z)h_1+a_2\sigm_-(z)h_2
+ a_{11}\Delta_{(11)}(\sigm_-(z)H_1)\\ &+&a_{21}\Delta_{(21)}(\sigm_-(z)H_1)
+\Delta_{(22)}(\sigm_-(z)H_1)+I_{1,4}\cr \cr
&=&a_0+a_1(x+z^{-1})+a_2(x^2+xz^{-1}+z^{-2})
\end{eqnarray*}
from which
\begin{eqnarray*}
\Gamma^\Vee_2(z)p(H_{2,4})&=&E_{1}(z)\sum_{\blamb\in\Pcal_{2,4}} a_\blamb\Delta_\blamb(\sigm_-(z)H_{1,4}))+I_{1,4}\\
&=&{1-xz\over z^2}\left(a_0+a_1\left(x+{1\over z}\right) +a_2\left(x^2+{x\over z}+{1\over z^2}\right) \right.+\\
&+&\left. a_{11}{x\over z}+a_{21}\left({x\over z^2}+{x^2\over
z}\right)+a_{22}{x^2\over z^2}\right).
\end{eqnarray*}
On the other hand
\begin{eqnarray*}
p(\ovsig_-(z) H_{3,4})&=&\sum_{\blamb\in \Pcal_{2,4}} a_\blamb\Delta_\blamb(\ovsig_-(z)H_{3,4})\\
&=&\sum_{\blamb\in\Pcal_{2,4}} a_{(\lambda_1,\lambda_2)}\det(h_{\lambda_j-j+i-1}-h_{\lambda_j-j+i}z^{-1})+I_{3,4}
\end{eqnarray*}
and so the equality
\be
\Gamma_2(z)p(H_{2,4})={z^2\over E_3(z)}\sum_{\blamb\in \Pcal_{2,4}} a_{(\lambda_1,\lambda_2)}\det(h_{\lambda_j-j+i-1}-h_{\lambda_j-j+i}z^{-1})+I_{3,4}\label{eq:beco}
\ee
holds in $B_{3,4}((z))$. Since
\[
{1\over E_3(z)}+I_{3,4}=1+h_1z+I_{3,4}=1+yz,
\]
formula~(\ref{eq:beco}) becomes:

\begin{eqnarray*}
\Gamma_2(z)p(H_{2,4})&=&z^2(1+yz)\left[a_0+a_1\left(y-{1\over z}\right)-a_2{y\over z}+a_{11}\left(y^2-{y\over z}+{1\over z^2}\right)\right.\\
&+& \left.
a_{21}\left({y\over z^2}-{y^2\over z}\right)+a_{22}{y^2\over z^2}\right].
\end{eqnarray*}
A few computations, carried out  by means of the {\tt CoCoA} software \cite{cocoa}, eventually tell
\begin{eqnarray*}
\Res_z(\Gamma_2(z)p(H_{2,4}))\cdot (\Gamma^\Vee_2(z)p(H_{2,4}))\hskip-7pt &=&\hskip-7pt(-a_{11} a_2 + a_1 a_{21}- a_0 a_{22})x^3
+ (a_{11} a_2  - a_1 a_{21} + a_0 a_{22}) x^2 y\\
& -&\hskip-6pt (a_{11} a_2  - a_1 a_{21} + a_0 a_{22}) x y^2+ (a_{11} a_2 - a_1 a_{21}  + a_0 a_{22})y^3\cr\cr
&=&(a_{11} a_2 - a_1 a_{21}  + a_0 a_{22})(y^3-y^2x+x^2y-x^3).
\end{eqnarray*}
So, the latter expression is identically zero if and only if the following  Pl\"ucker equation holds:
\be
a_{11} a_2 - a_1 a_{21}  + a_0 a_{22}=0.\label{eq:plucker}
\ee

\vspace{-8pt}
\section{The Grassmann Cone in  Infinite Exterior Power}\label{SecKP}

\begin{prop}\label{prop81} {\em  Let $r\geq 1$ be fixed. Then for all $(i_1,\ldots, i_r)\in \NN^r$, the equalities below hold in $B_r$:
\be
\hskip-1pt \sigm_-(z)(h_{i_1}\cdots h_{i_r})=\hskip-1pt\sigm_-(z)h_{i_1}\cdots \sigm_-(z)h_{i_r}\,\,\, \mathrm{and}\,\,\, \ovsig_-(z)(h_{i_1}\cdots h_{i_r})=\ovsig_-(z)h_{i_1}\cdots \ovsig_-(z)h_{i_r}.\label{eq8:qusring}
\ee
}
\end{prop}

\proof
Let us begin by checking  the first of~(\ref{eq8:qusring}).
Since $(\Delta_\blamb(H_r)\,|\, \blamb\in\Pcal_r)$ is a basis of $B_r$, each product of the form $h_{i_1}\cdots h_{i_r}$ is
 a linear combination $\sum_{\blamb\in\Pcal_{r,n}}a_\blamb\Delta_\blamb(H_r)$. Therefore
\[
\sigm_-(z)(h_{i_1}\cdots
h_{i_r})=\sigm_-(z)(\sum_{\blamb\in\Pcal_{r,n}}a_\blamb\Delta_\blamb(H_r))=\sum_{\blamb\in\Pcal_{r,n}}a_\blamb\sigm_-(z)\Delta_\blamb(H_r)=\sum_{\blamb\in\Pcal_{r,n}}a_\blamb\Delta_\blamb(\sigm_-(z)H_r),
\]
where last equality is due to Theorem~\ref{0stmnthm}. In other words:
\[
\sigm_-(z)(h_{i_1}\cdots h_{i_r})=\sum_{\blamb\in\Pcal_{r,n}}a_\blamb\det(\sigm_-(z)h_{\lambda_j-j+i})=\sigm_-(z)h_{i_1}\cdots \sigm_-(z)h_{i_r}.
\]

\vspace{-3pt}
\noindent
The proof for $\ovsig_-(z)$ is evidently analogous.\qed

\vspace{-6pt}
\claim{}\label{notasin} Let $B_\infty$ be the polynomial ring $\ZZ[e_1,e_2,\ldots]$ in infinitely many indeterminates $(e_1,e_2,\ldots)$.
 It is the projective limit of the rings $B_r$ in the category of graded $\ZZ$-algebras, in the following sense.
Recall notation and convention of Section~\ref{prelwtgrad}.  For all $s>r$ there are obvious projection maps $(B_s)_w\sra (B_r)_w$. 
They are  defined by $\Delta_\blamb(H_s)\mapsto\Delta_\blamb(H_r)$ if $\blamb\in\Pcal_r$ and by $\Delta_\blamb(H_s)\mapsto 0$ otherwise. Let
$
(B_\infty)_w:=\underset{\leftarrow}{\lim}(B_r)_w$. Clearly, for all $w\in \NN$ there exists $r>0$ such that $(B_\infty)_w=(B_r)_w$, the module of polynomials in $(e_1,\ldots,e_r)$ of weighted degree $w$ (See~\ref{prelwtgrad}). The
 ring $B_\infty$ is by definition the direct sum $\bigoplus_{w\geq 0}(B_\infty)_w$. In particular, for all $w\geq 0$ there exists $s\geq 0$ such that $(B_\infty)_{\leq w}:=\bigoplus_{0\leq i\leq w}(B_\infty)_i$ is isomorphic, as abelian group,  to $(B_r)_{\leq w}:=\bigoplus_{0\leq i\leq w}(B_r)_i$,  for all $r\geq s$.
 
Let $E_\infty(z)=1-e_1z+e_2z^2+\cdots$ and $\sum_{n\geq
0}h_nz^n=(E_\infty(z))^{-1}\in B_\infty[[z]]$. In this case the non-zero terms of the sequence
$H_\infty=(h_j)_{j\in\ZZ}$
 are algebraically independent and so $B_\infty$ is the free polynomial algebra $\ZZ[h_1,h_2,\ldots]$. Moreover $B_\infty=\bigoplus_{\blamb\in\Pcal_r}\Delta_\blamb(H_\infty)$.

\begin{corol}\label{corriso} {\em The maps  $\sigm_-(z),\ovsig_-(z):B_\infty\sra B_\infty[z^{-1}]$ are ring homomorphisms
 and are each other's inverses,  when regarded as elements of $(\End_\ZZ (B_\infty))[z^{-1}]$.
}
\end{corol}

\proof Consider an arbitrary product in $h_{i_1}\cdots h_{i_s}\in B_\infty$
 and let $w=i_1+\cdots+i_s$. There exists a sufficiently large $r>\max\{w, s\}$ such that $(B_\infty)_w=(B_r)_w$.
For such a choice of $r$, we have
\[
 \sigm_-(z)(h_{i_1}\cdots h_{i_s})= \sigm_-(z)h_{i_1}\cdots \sigm_-(z) h_{i_s}\qquad\mathrm{and}\qquad  \ovsig_-(z)(h_{i_1}\cdots h_{i_s})= \ovsig_-(z)h_{i_1}\cdots \ovsig_-(z) h_{i_s}
\]
by virtue of Proposition~\ref{prop81}.\qed

\vspace{-3pt}
\claim{\bf Proof of Corollary~\ref{mncor01}.}  Each $p\in B_\infty$ is a finite linear combination of $(\Delta_\blamb(H_\infty))_{\blamb\in \Pcal}$, where $\Pcal$ is the set of all the partitions. Then there exist $0\leq r_1\leq n$ such that  $p:=\sum_{\blamb\in\Pcal_{r_1,n}} a_\blamb\Delta_\blamb(H_\infty)$. Notice that the maximum weight of the partitions possibly occurring in the sum is $w:=r_1(n-r)$. Suppose that  there is $r_2\geq 0$ such that  $\phi_{r}(p)\in \Gcal_{r}$ for all $r\geq r_2$. Then  by Theorem~\ref{thm02}
\be
\Res_z \sum_{\blamb,\bmu\in \Pcal_{r_1,n}} a_\blamb a_\bmu E_{r-1}(z)\Delta_\blamb(\sigm_-(z)H_{r-1})\otimes {1\over E_{r+1}(z)}\ovsig_-(z)\Delta_\bmu(H_{r+1})=0, \label{eq8:rewr}
\ee
in $(B_{r-1})_{\leq w}\otimes (B_{r+1})_{\leq w}$ (notation as in~\ref{notasin}). Fix  $r\geq \max(r_1,r_2,w)$ big enough such that $(B_\infty)_{\leq w}$ is isomorphic to $(B_r)_{\leq w}$.  Then the left hand side of~(\ref{eq8:rewr}) is equal to
\be
\Res_z \sum_{\blamb,\bmu\in \Pcal} a_\blamb a_\bmu E_{r-1}(z)\Delta_\blamb(\sigm_-(z)H_{r-1})\otimes \big(\sum_{0\leq i\leq r }h_iz^i)\cdot \ovsig_-(z)\Delta_\bmu(H_{r+1}),\label{eq:xxxxxx}
\ee
because, to compute the residue, the formal power series  $\sum_{i\geq 0}h_iz^i=(E_{r+1}(z))^{-1}$ contributes with only finitely many summands. Therefore, invoking the isomorphism $(B_\infty)_{\leq w}\cong (B_r)_{\leq w}$, (\ref{eq:xxxxxx}) is equivalent to
\be
0=\Res_z \sum_{\blamb,\bmu\in \Pcal} a_\blamb a_\bmu E_{r-1}(z)\Delta_\blamb(\sigm_-(z)H_{\infty})\otimes \sum_{0\leq i\leq r }h_iz^i\cdot \ovsig_-(z)\Delta_\bmu(H_{\infty}), \label{eq8:rewrx}
\ee
in $(B_\infty\otimes B_\infty)((z))$.
 We may now replace $E_{r-1}(z)$ with $E_\infty(z)$ and $\sum_{0\leq i\leq r }h_iz^i$ with  $\sum_{i\geq 0}h_iz^i=(E_\infty(z))^{-1}$ in formula~(\ref{eq8:rewrx}), because adding powers of $z$ does not alterate the residue. We have so proven that  if $\phi_r(p)\in\Gcal_r$  for a big enough $r$,  then it satisfies equation~(\ref{eq1:mnthmm}).
Conversely, if $p$ satisfies (\ref{eq1:mnthmm}), then~(\ref{eq8:rewrx}) holds, and since the weight of the partitions involved to express $p$ as linear combination of Schur polynomials is bounded by a given positive integer $w$, there exists a big enough $s$ such that $(B_{r-1})_{\leq w}$ and $(B_{r+1})_{\leq w}$ are both isomorphic to $(B_\infty)_{\leq w}$, for all $r\geq s$. Thus equation~(\ref{eq8:rewrx}) is  equivalent to~(\ref{eq:xxxxxx}) and then to~(\ref{eq8:rewr}), i.e.  $\phi_r(p)\in \Gcal_r$ for all $r\geq s$.
As Corollary~\ref{corriso} implies the commutation $\Delta_\blamb(\sigm_-(z)H_\infty)=\sigm_-(z)\Delta_\blamb(H_\infty)$, for all $\blamb\in \Pcal$,  equation ~(\ref{eq8:rewr}) can be rewritten as 
\begin{eqnarray*}
\hskip70pt 0&=&\Res_z \sum_{\blamb\in \Pcal}   E_\infty(z)a_\blamb\sigm_-(z)\Delta_\blamb(H_\infty)\otimes {1\over E_\infty(z)}\sum_{\bmu\in \Pcal}a_\bmu\ovsig_-(z)\Delta_\bmu(H_\infty)\\
&=&\Res_z  E_\infty(z)\sigm_-(z)p\otimes {1\over E_\infty(z)}\ovsig_-(z)p.\hskip137pt \qed
\end{eqnarray*}


\claim{} Let $B:=B_\infty\otimes_\ZZ\QQ$ and define the sequence $\bfX:=(x_1,x_2,\ldots)$ through the equality
\[
\sum_{n\geq 0}h_nz^n=\exp(\sum_{i\geq 1}x_iz^i),
\]
holding in $B[[z]]$, in such a way that  each $h_n$ can be regarded as a function of $(x_1,x_2,\ldots)$. Standard calculations show that $h_n$
 is a polynomial expression of $(x_1,\ldots, x_n)$, homogeneous of degree $n$ with respect to the weight graduation of $B_\infty$ (for which $h_n$ and $x_n$ have degree $n$). 
\bclm{\bf Lemma.} {\em The following equalities hold in the ring  $B$  for all $j\geq 1$ and $n\geq 0$:
\be
{\partial^j h_n\over \d x_1^j}={\partial h_n\over\d x_j}=h_{n-j}.\label{eq8:lmhn-j}
\ee
}
\eclm
\proof
For all $j\geq 1$:
\[
\sum_{n\geq 0}{\d h_n \over \d  x_j}z^n={\d \over \d x_j}\sum_{n\geq 0}h_nz^n={\d \over \d x_j}\exp(\sum_{j\geq1}x_jz^j) =z^j\exp(\sum_{j\geq1}x_jz^j)=\sum_{n\geq 0}h_nz^{n+j}.
\]
Comparing the coefficient of $z^n$ in the first and last side gives
\[
{\d h_n \over \d x_j}=h_{n-j}.
\]
In particular  $\d h_n/\d x_1=h_{n-1}$.  Iterating  $j$ times the operator  $\d/\d x_1$ gives~(\ref{eq8:lmhn-j}), as desired.\qed

Let
\[
\Gamma_\infty(z)\Delta_\blamb(H_\infty):={1\over E_\infty(z)}\ovsig_-(z)\Delta_\blamb(H_\infty)\in B_\infty((z))
\]
and
\[
\Gamma_\infty^\Vee(z)\Delta_\blamb(H_\infty):=E_\infty(z)\sigm_-(z)\Delta_\blamb(H_\infty)\in B_\infty((z)).
\]

\medskip
Define  $\Gamma(z), \Gamma^\Vee(z)$ to be,  respectively,  $\Gamma_\infty(z)\otimes {1_{\QQ}}:B\sra B((z))$ and $\Gamma^\Vee_\infty(z)\otimes {1_{\QQ}}:B\sra B((z))$.
Corollary~\ref{mncor02}  follows immediately from Corollary~\ref{mncor01} and the following:
\begin{thm}[Cf. \cite{KR}, Theorem 5.1] {\em We have: \be \Gamma(z)=\exp\left(\sum_{i\geq
1}x_iz^i\right)\exp\left(-\sum_{i\geq 1}{1\over iz^i}{\d\over \d
x_i}\right)\label{eq8:vop} \ee and \be
\Gamma^\Vee(z)=\exp\left(-\sum_{i\geq
1}x_iz^i\right)\exp\left(\sum_{i\geq 1}{1\over iz^i}{\d\over \d
x_i}\right).\label{eq8:vopv} \ee} \end{thm}

\proof
Since
\[
{1\over E_\infty(z)}=\sum_{n\geq 0}h_nz^n=\exp(\sum_{i\geq 1}x_iz^i),
\]
it follows that $E_\infty(z)=\exp(-\sum_{i\geq 1}x_iz^i)$  and   the first factors involved on the left hand side of~(\ref{eq8:vop}) and~(\ref{eq8:vopv})
 are explained. Let us now observe that:
\[
\ovsig_-(z)h_n=h_n-{h_{n-1}\over z}=\left(1-{1\over z}{\d \over \d x_1}\right)h_n .
\]
Evaluating the  well--known identity  $1-t=\displaystyle{\exp\left(-\sum_{n\geq 1}{t^n\over n}\right)}$ at $t=z^{-1}\displaystyle{\d\over \d x_1}$, and using~(\ref{eq8:lmhn-j}), we have
\be
\ovsig_-(z)h_n=\exp\left(-\sum_{i\geq 1}{1\over iz^i}{\d^i\over \d x_1^i}\right)h_n=\exp\left(-\sum_{i\geq 1}{1\over iz^i}{\d\over \d x_i}\right)h_n.\label{eq:yyy}
\ee
Now, we observe that
$
\displaystyle{\exp\left(-\sum_{i\geq 1}{1\over iz^i}{\d\over \d x_i}\right):B\sra B[z^{-1}]}
$
is a ring homomorphism, because it is  the exponential of the  first order differential operator
$
\displaystyle{-\sum_{i\geq 1}{1\over i z^i}{\d\over \d x_i}.}
$
Thus
\[
\ovsig_-(z)=\exp\left(-\sum_{i\geq 1}{1\over iz^i}{\d\over \d
x_i}\right),
\]
because  both sides are ring homomorphisms and by~(\ref{eq:yyy}) they coincide on $h_n$, for all $n\geq 0$, which generate $B$ as a $\QQ$-algebra. The proof  that
\[
\sigm_-(z)=\exp\left(\sum_{i\geq 1}{1\over iz^i}{\d\over \d x_i}\right)
\]
is similar,  but  arguing  that $\sigm_-(z)$ is the inverse of $\ovsig_-(z)$ in $\End_\QQ(B)[z^{-1}]$ turns it easier. \qed

\begin{rmk} As pointed out in the introduction, for each $\tau$-function, the normalized first order formal pseudo-differential operator $L_\tau:=P_\tau(\d)\d P_\tau(\d)^{-1}$, where
\[
P_\tau(z)={\ovsig_-(z)\tau\over \tau}\in B_{(0)}[z^{-1}]
\]
is a solution of the KP hierarchy in Lax form~(\ref{eq:lax}). See e.g.~\cite[Section 7.5]{KR}.
\end{rmk}

\claim{}\label{weight} Let $(\bw^rM_0)_w:=\bigoplus_{|\blamb|=w}\ZZ\cdot\wb^r_\blamb$. The isomorphism $B_r\sra \bw^rM_0$ does restrict to an isomorphism $(B_r)_w\sra (\bw^rM_0)_w$ and the map $(B_{r_1})_w\sra (B_{r_2})_w$, for all $r_1\geq r_2$, induce the  epimorphism $(\bw^{r_1}M_0)_w\sra (\bw^{r_2}M_0)_w$ mapping $\wb^{r_1}_\blamb\mapsto \wb^{r_2}_\blamb$ (where by convention $\wb^r_\blamb=0$ if $\ell(\blamb)\geq r$).  The projective limit  $(\bw^\infty M_0)_w:=\underset{\leftarrow}{\lim}(\bw^rM_0)_w$, with respect to the above projection maps,  may be identified with the free abelian group generated by the symbols $\wb^\infty_\blamb$, where $\blamb$ ranges over all the partitions of weight $w$. Let $\bw^\infty M_0=\bigoplus_{w\geq 0}(\bw^\infty M_0)_w=\bigoplus_{\blamb\in\Pcal}\ZZ\wb^\infty_\blamb$ and
\vspace{-2pt}
$$
\phi_\infty: B_\infty\sra \bw^\infty M_0
$$
be the $\ZZ$-isomorphism defined by $\Delta_\blamb(H_\infty)\mapsto \wb^\infty_\blamb$ (the boson-fermion correspondence).   
Then  Corollaries~\ref{mncor01} and~\ref{mncor02} show that one may safely define 
 the {\em Grassmann cone} $\Gcal_\infty\subseteq \bw^\infty M_0$ as the locus of all $\bfm\in\bigoplus_{\blamb\in\Pcal}\ZZ\wb^\infty_\blamb$ such that $\phi_\infty^{-1}(\bfm)$ satisfies~(\ref{eq1:mnthmm}). Similarly, the Grassmann cone $\Gcal_\infty\otimes \QQ$ is the locus of $\bfm\in \bw^\infty(M_0\otimes\QQ)$ such that $(\phi_\infty\otimes 1)^{-1}(\bfm)\in B$ is a {\em tau}  function for the $KP$ hierarchy~(\ref{eq:KPHwVOx}).

\medskip
\newpage


\parbox[t]{3in}{{\rm Letterio~Gatto}\\
{\tt \href{mailto:letterio.gatto@polito.it}{letterio.gatto@polito.it}}\\
{\it Dipartimento~di~Scienze~Matematiche}\\
{\it Politecnico di Torino}\\
{\it ITALY}} \hspace{1.5cm}
\parbox[t]{2.5in}{{\rm Parham~Salehyan}\\
{\tt \href{mailto:parham@ibilce.unesp.br}{parham@ibilce.unesp.br}}\\
{\it Ibilce UNESP}\\
{\it Campus de S\~ao Jos\'e do Rio Preto, SP}\\
{\it BRAZIL}}

\end{document}